\documentclass[a4paper,12pt]{amsart}
\usepackage{amssymb,amscd,amsmath,a4wide,graphicx,stmaryrd,fullpage,setspace,microtype,xr,tikz,textcomp,enumitem,wasysym,cancel}
\usepackage{hyperref}
\definecolor{dark-red}{rgb}{0.5,0.15,0.15}
\hypersetup{
	colorlinks   = true, 
	urlcolor     = dark-red, 
	linkcolor    = black, 
	citecolor   = dark-red 
}

\usepackage[all]{xy}
\usepackage[T1]{fontenc}
\usepackage{lmodern}
\usepackage[numbers,sort&compress]{natbib}
\title{Regular directed path and Moore flow}

\author[P. Gaucher]{Philippe Gaucher}

\address{Universit\'e Paris Cit\'e, CNRS, IRIF, F-75013, Paris, France}

\urladdr{http://www.irif.fr/{\~{}}gaucher} 

\makeatletter
\@namedef{subjclassname@2020}{\textup{2020} Mathematics Subject Classification}
\makeatother


\swapnumbers

\newcommand{\C}{\mathcal{C}}
\newcommand{\D}{\mathcal{D}}
\newcommand{\K}{\mathcal{K}}

\newcommand{\W}{\mathcal{W}}
\newcommand{\F}{\mathcal{F}}

\newcommand{\de}{\partial}
\newcommand{\p}{\times}
\renewcommand{\vec}{\overrightarrow}
\renewcommand{\P}{\mathbb{P}}

\newtheorem*{thmN}{Theorem}

\newtheorem{thm}{Theorem}[section]
\newtheorem{prop}[thm]{Proposition}
\newtheorem{lem}[thm]{Lemma}
\newtheorem{conj}[thm]{Conjecture}
\newtheorem{cor}[thm]{Corollary}

\newtheorem{defnot}[thm]{Definition and notation}
\newcommand{\bdn}{\begin{defnot}}
	\newcommand{\edn}{\end{defnot}}
\newcommand{\bp}{\begin{prop}}
	\newcommand{\ep}{\end{prop}}
\newcommand{\bth}{\begin{thm}}
	\renewcommand{\eth}{\end{thm}}
\newcommand{\bpf}{\begin{proof}}
	\newcommand{\epf}{\end{proof}}
\newcommand{\bc}{\begin{cor}}
	\newcommand{\ec}{\end{cor}}

\theoremstyle{definition}
\newtheorem{defn}[thm]{Definition}
\newtheorem{rem}[thm]{Remark}
\newtheorem{exa}[thm]{Example}
\newcommand{\bd}{\begin{defn}}
	\newcommand{\ed}{\end{defn}}
\newtheorem{nota}[thm]{Notation}

\newcommand{\fL}[1]{\ar@{->}[ll]_-{#1}}
\newcommand{\fR}[1]{\ar@{->}[rr]^-{#1}}
\newcommand{\fRr}[1]{\ar@{->}[rrr]^-{#1}}
\newcommand{\fD}[1]{\ar@{->}[dd]_-{#1}}
\newcommand{\fU}[1]{\ar@{->}[uu]^-{#1}}
\newcommand{\f}[2]{\ar@{->}[#1]|{#2}}
\newcommand{\ff}[2]{\ar@2{->}[#1]|{#2}}
\newcommand{\frr}[1]{\ar@{->}[rrrr]^-{#1}}

\newcommand{\fl}[1]{\ar@{->}[l]_-{#1}}
\newcommand{\fr}[1]{\ar@{->}[r]^-{#1}}
\newcommand{\fd}[1]{\ar@{->}[d]_-{#1}}
\newcommand{\fu}[1]{\ar@{->}[u]^-{#1}}

\renewcommand{\top}{{\mathbf{Top}}}

\newcommand{\iso}{\cong}

\newcommand{\ot}{\otimes}

\newcommand{\vI}{\vec{I}}
\renewcommand{\leq}{\leqslant}
\renewcommand{\geq}{\geqslant}

\newcommand{\ptop}[1]{{\brm{{#1}dTop}}}

\newcommand{\moore}{{\mathbb{M}}}
\newcommand{\lmoore}{\mathbb{M}_!}

\def\cartesien{%
	\ar@{-}[]+R+<6pt,-2pt>;[]+RD+<6pt,-6pt>%
	\ar@{-}[]+D+<2pt,-6pt>;[]+RD+<6pt,-6pt>%
}
\def\cocartesien{%
	\ar@{-}[]+L+<-6pt,+2pt>;[]+LU+<-6pt,+6pt>%
	\ar@{-}[]+U+<-2pt,+6pt>;[]+LU+<-6pt,+6pt>%
}
\def\hocartesien{%
	\ar@{-}[]+R+<6pt,-2pt>;[]+RD+<6pt,-6pt>_{h}%
	\ar@{-}[]+D+<2pt,-6pt>;[]+RD+<6pt,-6pt>%
}
\def\hococartesien{%
	\ar@{-}[]+L+<-6pt,+2pt>;[]+LU+<-6pt,+6pt>_{h}%
	\ar@{-}[]+U+<-2pt,+6pt>;[]+LU+<-6pt,+6pt>%
}

\newcommand{\brm}[1]{{\rm{\mathbf{#1}}}}

\newcommand{\dtop}{{\brm{Flow}}}

\newcommand{\dtopG}{{\mathcal{G}\brm{Flow}}}

\newcommand{\set}{{\brm{Set}}}

\newcommand{\ttop}{{\brm{TOP}}}

\newcommand{\glob}{{\mathrm{Glob}}}
\newcommand{\topdgr}{[\mathcal{G}^{op},\top]}
\newcommand{\topdgrq}{[\mathcal{G}^{op},\top_q]_0^{proj}}
\newcommand{\topdgrm}{[\mathcal{G}^{op},\top_m]_0^{proj}}
\newcommand{\topdgrh}{[\mathcal{G}^{op},\top_h]_0^{proj}}

\newcommand{\ttt}{two-out-of-three property}

\DeclareMathOperator{\id}{Id}

\DeclareMathOperator{\Obj}{Obj}
\DeclareMathOperator{\Mor}{Mor}

\DeclareMathOperator{\Ch}{Ch}
\DeclareMathOperator{\Cat}{\mathrm{cat}}
\newcommand{\mins}[1]{\texorpdfstring{$#1$}{Lg}}

\newcommand{\liminj}{\varinjlim}
\newcommand{\limproj}{\varprojlim}

\newcommand{\rest}{\!\upharpoonright\!}

\renewcommand{\P}{\mathbb{P}}

\makeatletter
\def\varholim@#1#2{%
	\vtop{\m@th\ialign{##\cr
			\hfil$#1\operator@font holim$\hfil\cr
			\noalign{\nointerlineskip\kern1.5\ex@}#2\cr
			\noalign{\nointerlineskip\kern-\ex@}\cr}}%
}
\def\holimproj{%
	\mathop{\mathpalette\varholim@{\leftarrowfill@\textstyle}}\nmlimits@
}
\def\holiminj{%
	\mathop{\mathpalette\varholim@{\rightarrowfill@\textstyle}}\nmlimits@
}
\makeatother

\newcommand{\ddownarrow}{{\downarrow}}

\DeclareMathOperator{\pf}{{Gph}}
\DeclareMathOperator{\cocyl}{{Path}}
\newcommand{\adj}[5]{\xymatrix@C=#5em{{#1}\ar@/^0.8em/[r]^-{#2} \ar@{}[r]|-{\perp} & \ar@/^0.8em/[l]^-{#3} {#4}}}

\DeclareMathOperator{\seq}{Seq}

\allowdisplaybreaks

\begin{document}

\begin{abstract}
	Using the notion of tame regular $d$-path of the topological $n$-cube, we introduce the tame regular realization of a precubical set as a multipointed $d$-space. Its execution paths correspond to the nonconstant tame regular $d$-paths in the geometric realization of the precubical set. The associated Moore flow gives rise to a functor from precubical sets to Moore flows which is weakly equivalent in the h-model structure to a colimit-preserving functor. The two functors coincide when the precubical set is spatial, and in particular proper. As a consequence, it is given a model category interpretation of the known fact that the space of tame regular $d$-paths of a precubical set is homotopy equivalent to a CW-complex. We conclude by introducing the regular realization of a precubical set as a multipointed $d$-space and with some observations about the homotopical properties of tameness.
\end{abstract}
	
\maketitle
\tableofcontents
\hypersetup{linkcolor = dark-red}

\section{Introduction}

\subsection*{Presentation}
It is described in \cite{NaturalRealization} a way of realizing a precubical set as a flow without any non-canonical choice of any cofibrant replacement in the construction, by introducing a natural realization functor. It is an improvement of what is done in \cite[Section~4]{realization}. The latter terminology comes from the fact that it uses the notion of \textit{natural $d$-path} of the topological $n$-cube. This work is devoted to the study of a way of realizing a precubical set as a multipointed $d$-space without any non-canonical choice in the construction either, unlike what is done in \cite[Section~5]{realization}. Precubical sets are a model for concurrency theory \cite{DAT_book}. The $n$-cube represents the concurrent execution of $n$ actions. This idea is even further developed in \cite{zbMATH07226006} in which it is established that the non-positively curved precubical sets contain most of the useful examples. However, we do not need this concept and we keep working with general precubical sets like in \cite{NaturalRealization}. The $d$-paths of a precubical set, i.e. in the geometric realization of the precubical set, are the continuous paths which are compositions of continuous paths included in cubes of the geometric realization which are locally nondecreasing with respect to each axis of coordinates. These particular continuous paths represent the possible execution paths of the concurrent process modelled by the precubical set. The local nondecreasingness represents time irreversibility. The notion of regular continuous path is introduced in \cite[Definition~1.1]{reparam}. The terminology comes from elementary differential geometry where a differentiable maps $\gamma:[0,1]\to \mathbb{R}^n$ is called regular if $p'(t)\neq 0$ for all $t\in ]0,1[$. Intuitively, it is a continuous path without zero speed point.

By starting from the topological $n$-cube equipped with the set of nonconstant tame regular $d$-paths of the topological $n$-cube, it is introduced the \textit{tame regular realization} $|K|^t_{reg}$ of a precubical set $K$ as a multipointed $d$-space in Definition~\ref{treg_rea}. Since the composition of two tame regular $d$-paths is still tame regular, all execution paths of this multipointed $d$-space are tame regular $d$-paths (and nonconstant). Thus, the tame regular realization of a precubical set does not contain as execution paths all nonconstant regular $d$-paths between two vertices of the precubical set: see Figure~\ref{tame-nontame} for an example of a non-tame $d$-path. This point of view is not restrictive: see Theorem~\ref{comparison-tame-nottame-Moore-flow}, Corollary~\ref{comparison-tame-nottame-Moore-flow-dspace} and Corollary~\ref{reg_weak}.

It is well established that the relevant information for concurrency theory is contained in the homotopy type of the space of execution paths \cite{DAT_book}. In particular, it is not contained in the topology of the underlying space. It therefore suffices to consider the Moore flow $\moore^{\mathcal{G}}(|K|^t_{reg})$. The latter is obtained by forgetting the underlying topological space of the multipointed $d$-space $|K|^t_{reg}$ and by keeping only the execution paths and the information about their reparametrization (see Theorem~\ref{MG}). 

The point is that the functor $K\mapsto\moore^{\mathcal{G}}(|K|^t_{reg})$ is the composite of a left adjoint $K\mapsto |K|^t_{reg}$ from precubical sets to multipointed $d$-spaces and of a right adjoint $\moore^{\mathcal{G}}:\ptop{\mathcal{G}} \to \dtopG$ from multipointed $d$-spaces to Moore flows. The latter does not commute with colimits in general. The main result of this paper is that the composite functor $K\mapsto \moore^{\mathcal{G}}(|K|^t_{reg)}$ from precubical sets to Moore flows can be replaced, up to a natural weak equivalence of the h-model structure of Moore flows, by a colimit-preserving functor. This replacement up to a weak equivalence of the h-model structure means that the spaces of execution paths are replaced by homotopy equivalent ones. Moreover, the image of this new functor is included in the class of m-cofibrant Moore flows. More precisely, we obtain the two following theorems.

\begin{thmN}  (Theorem~\ref{iso_reg_reg0})
	There exists a colimit-preserving functor \[[-]_{reg}:\square^{op}\set \longrightarrow \dtopG\] from precubical sets to Moore flows and a natural map of Moore flows \[[K]_{reg} \longrightarrow \moore^{\mathcal{G}}(|K|^t_{reg})\] which is a weak equivalence of the h-model structure of Moore flows for all precubical sets $K$. Moreover, the above natural map of Moore flows is an isomorphism if and only if $K$ is spatial. 
\end{thmN}

\begin{thmN} (Theorem~\ref{m-cof-moore-flow})
	For all precubical sets $K$, the Moore flow $[K]_{reg}$ is m-cofibrant.
\end{thmN}

As examples of spatial precubical sets, there are all proper precubical sets by \cite[Proposition~7.5]{NaturalRealization}, and in particular the geometric and non-positively curved ones. In other terms, the precubical sets coming from a lot of real concurrent systems by \cite[Proposition~1.29]{zbMATH07226006}. 

As any right adjoint between locally presentable category, the functor $\moore^{\mathcal{G}}:\ptop{\mathcal{G}} \to \dtopG$ is accessible. However it is not finitely accessible. This is due to the fact that the convenient category of topological spaces we are working with is locally $(2^{\aleph_0})^+$-presentable and not locally finitely presentable. However, it does commute with some particular transfinite compositions \cite[Theorem~6.14]{Moore1} \cite[Theorem~5.7]{Moore2}. The right adjoint $\moore^{\mathcal{G}}:\ptop{\mathcal{G}} \to \dtopG$ is also a part of a Quillen equivalence $\lmoore^{\mathcal{G}}\dashv \moore^{\mathcal{G}}$ between the q-model structures of multipointed $d$-spaces and Moore flows by \cite[Theorem~8.1]{Moore2} which satisfies a very specific property: the unit and the counit of this Quillen equivalence induce isomorphisms between the q-cofibrant objects \cite[Theorem~7.6 and Corollary~7.9]{Moore2}. These facts together with the isomorphism of Theorem~\ref{iso_reg_reg0} when $K$ is a spatial precubical set is an illustration of the following informal observation: the right adjoint $\moore^{\mathcal{G}}:\ptop{\mathcal{G}} \to \dtopG$ commutes with good enough colimits.

It is then deduced from Theorem~\ref{m-cof-moore-flow} that the space of tame regular $d$-paths of a precubical set is homotopy equivalent to a CW-complex. We obtain a purely model category proof of this known fact (see the long comment after Corollary~\ref{space-mcof} for bibliographical references).

\begin{thmN} (Corollary~\ref{space-mcof})
For every precubical set $K$, the space of tame regular $d$-paths is homotopy equivalent to a CW-complex.
\end{thmN}

At this point, the reader may wonder what happens to the theory developed in this paper when the tameness condition is removed. As an answer, we introduce the \textit{regular realization} $|K|_{reg}$ of a precubical set $K$ as a multipointed $d$-space in Definition~\ref{reg_rea}. The execution paths of $|K|_{reg}$ are exactly all nonconstant regular $d$-paths in the geometric realization of $K$ between two vertices of $K$, not only the tame ones. We then obtain the following comparison theorem as a corollary of Theorem~\ref{iso_reg_reg0} and of results from Raussen and Ziemia{\'{n}}ski.

\begin{thmN} (Corollary~\ref{reg_weak})
	For all precubical sets $K$, there exists a natural weak equivalence of the h-model structure of Moore flows 
	\[
	[K]_{reg} \longrightarrow \moore^{\mathcal{G}}(|K|_{reg}).
	\]
\end{thmN}

The notion of natural $d$-path of a precubical set enables us in \cite{NaturalRealization} to find a way to realize any precubical set as an m-cofibrant flow without using any non-canonical choice of any cofibrant replacement functor on the category of flows. In this paper, the notion of tame regular $d$-path of a precubical set enables us to obtain a realization of any precubical set as a m-cofibrant Moore flow without using any non-canonical choice of any cofibrant replacement functor on the category of multipointed $d$-spaces or on the category of Moore flows either. 

In both cases, we obtain much better results than in \cite{realization} because all constructions become canonical thanks to the geometric properties of the cubes. However, the associated Moore flow is m-cofibrant, and not q-cofibrant. It is the counterpart of these improvements: unlike in \cite{realization}, we have to deal with accessible model categories which are unlikely to be combinatorial. 

By removing the adjective \textit{regular} from Definition~\ref{treg_rea}, it is possible to develop the same theory as the one develop in this paper. However, it is necessary to replace the reparametrization category $\mathcal{G}$ defined in Notation~\ref{defGcat} by the reparametrization category $\mathcal{M} \supset \mathcal{G}$ of \cite[Proposition~4.11]{Moore1}: the small category $\mathcal{G}$ is exactly the subcategory of $\mathcal{M}$ with the same set of objects $]0,+\infty[$ and containing only the invertible maps of $\mathcal{M}$. In particular, it is necessary to work with $\mathcal{M}$-flows instead of with $\mathcal{G}$-flows (the latter are called Moore flows in this paper). Some results in this direction are already available in \cite{Moore3}. 

\subsection*{Outline of the paper}

Section~\ref{reg_sec}, after a reminder about multipointed $d$-spaces and precubical sets, introduces the \textit{tame regular realization} of a precubical set as a multipointed $d$-space. It is built using the notion of \textit{nonconstant tame regular $d$-path} of the topological $n$-cube. The execution paths of the tame regular realization are exactly the nonconstant tame regular $d$-paths in the geometric realization of the precubical set.

Section~\ref{qhm_sec} starts by a reminder about \textit{flows} and \textit{Moore flows} (flows are just a particular case of Moore flows for which the reparametrization category is the terminal category), as well as some results about them which are used in this paper. We add in this reminder the construction of the $\{q,m,h\}$-model structures of Moore flows. It is an easy adaptation of the case of flows treated in \cite[Theorem~7.4]{QHMmodel}.

Section~\ref{length_sec} is devoted to the important notion of \textit{$L_1$-arc length} of a $d$-path in the geometric realization of a precubical set. Some results coming from \cite{MR2521708} are adapted to $\Delta$-generated spaces. Theorem~\ref{Psi} is a slightly improved statement of a similar statement in \cite{MR2521708}: a homotopy equivalence is replaced by a homeomorphism. It is then proved that Theorem~\ref{Psi} implies that the flow associated to the tame regular realization as a multipointed $d$-space is the \textit{tame concrete realization} of a precubical set as recalled in Definition~\ref{rea_conc} and originally defined in \cite[Definition~7.1]{NaturalRealization}. 

Section~\ref{cubechain} makes the link between the ideas of this paper and Ziemia{\'{n}}ski notion of cube chain. It is an adaptation of \cite[Section~6]{NaturalRealization} in the setting of Moore flows. The idea is to adapt Ziemia{\'{n}}ski's cube chains initially developed for the closed monoidal category of topological spaces $(\top,\p)$ to the biclosed semimonoidal category of $\mathcal{G}$-spaces $(\topdgr_0,\ot)$. It culminates with Theorem~\ref{doublebarG}. 

Section~\ref{main}, after recalling the important notion of \textit{spatial precubical set}, establishes the main results of this paper, namely Theorem~\ref{iso_reg_reg0}, Theorem~\ref{m-cof-moore-flow} and Corollary~\ref{space-mcof} expounded in the introduction. 

Finally, Section~\ref{question-r} addresses the question of the \textit{regular realization} of a precubical set as a multipointed $d$-space and a comparison theorem with the tame regular realization is obtained in Theorem~\ref{comparison-tame-nottame-Moore-flow} and Corollary~\ref{comparison-tame-nottame-Moore-flow-dspace}. We then deduce Corollary~\ref{reg_weak} from the results of this section and Theorem~\ref{iso_reg_reg0}. Some observations are also made about the link between tameness and model categories. 


\subsection*{Prerequisites and notations}

We refer to \cite{TheBook} for locally presentable categories, to \cite{MR2506258} for combinatorial model categories.  We refer to \cite{MR99h:55031} and to \cite{ref_model2} for more general model categories, and to \cite{zbMATH06722019,HKRS17,GKR18} for accessible model categories. We refer to \cite{KellyEnriched} and to \cite[Chapter~6]{Borceux2} for enriched categories. All enriched categories are topologically enriched categories: \textit{the word topologically is therefore omitted}. 

The initial object of a category is denoted by $\varnothing$. The terminal object of a category is denoted by $\mathbf{1}$. The set of maps from $X$ to $Y$ in a category $\C$ is denoted by $\C(X,Y)$. $\id_X$ denotes the identity map of $X$. Let $\mathcal{I}$ be a set of maps of a bicomplete category $\C$. An \textit{$\mathcal{I}$-cellular} object is an object $X$ such that the canonical map $\varnothing \to X$ is a transfinite composition of pushouts of maps of $\mathcal{I}$. This sequence of pushouts is called the \textit{$\mathcal{I}$-cellular decomposition} of $X$. An \textit{$\mathcal{I}$-cofibrant} object is a retract of an $\mathcal{I}$-cellular object. A \textit{cellular} (\textit{cofibrant} resp.) object of a combinatorial model category is an $\mathcal{I}$-cellular ($\mathcal{I}$-cofibrant resp.) object where $\mathcal{I}$ is the set of generating cofibrations. The prefix $\mathcal{I}$ can be omitted when there is no ambiguity. 

The category $\top$ denotes the category of \textit{$\Delta$-generated spaces} or of \textit{$\Delta$-Hausdorff $\Delta$-generated spaces} together with the continuous maps (cf. \cite[Section~2 and Appendix~B]{leftproperflow}). The inclusion functor from the full subcategory of $\Delta$-generated spaces to the category of general topological spaces together with the continuous maps has a right adjoint called the $\Delta$-kelleyfication functor and denoted by $k_\Delta$. The latter functor does not change the underlying set. It only adds open subsets. The category $\top$ is locally presentable by \cite[Theorem~3.9]{FR} and Cartesian closed by \cite[Theorem~3.6]{MR45:9323} (see also \cite[Proposition~2.5]{mdtop}). The internal hom $\ttop(X,Y)$ is given by taking the $\Delta$-kelleyfication of the compact-open topology on the set $\top(X,Y)$. The category $\top$ can be equipped with its q-model structure $\top_q$. The m-model structure $\top_m$ \cite{mixed-cole} and the h-model structure $\top_h$ \cite{Barthel-Riel} of $\top$ are also used in various places of the paper. Compact means quasicompact Hausdorff (French convention). 

The following elementary fact is implicitly assumed several times in this work to use results written with the compact-open topology and not with its $\Delta$-kelleyfication: 

\bp \label{homotopy-deltakelleyfication}
Let $X$ and $Y$ be two homotopy equivalent general topological spaces. Then their $\Delta$-kelleyfications $k_\Delta(X)$ and $k_\Delta(Y)$ are homotopy equivalent.
\ep

\bpf
It suffices to prove that two homotopic maps $f,g:X\rightrightarrows Y$ give rise to two homotopic maps $k_\Delta(f),k_\Delta(g):k_\Delta(X)\rightrightarrows k_\Delta(Y)$. Let $H:X\p [0,1]\to Y$ be a homotopy from $f$ to $g$. Since $k_\Delta$ is a right adjoint, it preserves binary products. Therefore there are the natural homeomorphisms $k_\Delta(X\p [0,1]) \iso k_\Delta(X) \p k_\Delta([0,1]) \iso k_\Delta(X)\p [0,1]$, the right-hand one because $[0,1]$ is $\Delta$-generated. Thus $k_\Delta(H)$ yields a homotopy from $k_\Delta(f)$ to $k_\Delta(g)$.
\epf 

We conclude the prerequisites by giving the general definition of a \textit{reparametrization category} for the convenience of the reader because the terminology is used sometimes in this paper. It is not required to understand the proofs of this paper.

\bd A \textit{semimonoidal category} $(\K,\ot)$ is a category $\K$ equipped with a functor $\ot:\K\p \K\to \K$ together with a natural isomorphism $a_{x,y,z}: (x\ot y) \ot z \to x \ot (y\ot z)$ called the \textit{associator} satisfying the pentagon axiom \cite[diagram~(5) page 158]{MR1712872}. 
\ed

According to the usual terminology used for similar situations, a semimonoidal category could be called a \textit{non-unital monoidal category}. Note that it is the monoidal structure which is non unital, not the category.

\bd A semimonoidal category $(\K,\ot)$ is \textit{enriched} if the category $\K$ is enriched and if the set map \[\K(a,b)\p \K(c,d) \longrightarrow \K(a\ot c,b\ot d)\] is continuous for all objects $a,b,c,d\in \Obj(\K)$. 
\ed

\bd \label{def-reparam} \cite[Definition~4.3]{Moore1}
A \textit{reparametrization category} $(\mathcal{P},\ot)$ is a small enriched semimonoidal category satisfying the following additional properties: 
\begin{enumerate}
	\item The semimonoidal structure is strict, i.e. the associator is the identity.
	\item All spaces of maps $\mathcal{P}(\ell,\ell')$ for all objects $\ell$ and $\ell'$ of $\mathcal{P}$ are contractible. 
	\item For all maps $\phi:\ell\to \ell'$ of $\mathcal{P}$, for all $\ell'_1,\ell'_2\in \Obj(\mathcal{P})$ such that $\ell'_1\ot\ell'_2=\ell'$, there exist two maps $\phi_1:\ell_1\to \ell'_1$ and $\phi_2:\ell_2\to \ell'_2$ of $\mathcal{P}$ such that $\phi=\phi_1 \ot \phi_2 : \ell_1\ot\ell_2 \to \ell'_1 \ot\ell'_2$ (which implies that $\ell_1 \ot \ell_2=\ell$). 
\end{enumerate}
\ed 

\begin{nota}
	To stick to the intuition, we set $\ell+\ell' := \ell \ot \ell'$ for all $\ell,\ell'\in \Obj(\mathcal{P})$. Indeed, morally speaking, $\ell$ is the length of a path.
\end{nota}

A reparametrization category $\mathcal{P}$ is an enriched category with contractible spaces of morphisms such that  the set $\Obj(\mathcal{P})$ of objects of $\mathcal{P}$ has a structure of a semigroup with a composition law denoted by $+$, such that the set map 
\[
\ot : \mathcal{P}(\ell_1,\ell'_1) \p \mathcal{P}(\ell_2,\ell'_2) \to \mathcal{P}(\ell_1+\ell_2,\ell'_1+\ell'_2)
\]
is continuous for all $\ell_1,\ell'_1,\ell_2,\ell'_2 \in \Obj(\mathcal{P})$, and such that every map of $\mathcal{P}$ is of the form $\phi_1 \ot \phi_2$ (not necessarily in a unique way).

\begin{exa} \label{terminal-reparam}
	The terminal category with one object $\underline{1}$ and one map $\id_{\underline{1}}$ is a reparametri\-zation category. 
\end{exa}

\subsection*{Acknowledgment}

I would like to thank the anonymous referee for all comments, which enabled me to substantially improve the paper.

\section{The tame regular realization of a precubical set}
\label{reg_sec}

This section starts by a reminder about multipointed $d$-spaces, then another one about precubical sets. Then we introduce the notion of $d$-path and, finally, the tame regular realization of a precubical set is defined.

\begin{nota}
	The notations $\ell,\ell',\ell_i,L,\dots$ mean a strictly positive real number unless specified something else.
\end{nota}

\begin{nota}
	Let $\ell>0$. Let $\mu_{\ell}:[0,\ell]\to [0,1]$ be the homeomorphism defined by $\mu_\ell(t) = t/\ell$. 
\end{nota}

\begin{nota} \label{defG}
	The notation $[0,\ell_1]\iso^+ [0,\ell_2]$ means a nondecreasing homeomorphism from $[0,\ell_1]$ to $[0,\ell_2]$ with $\ell_1,\ell_2>0$. Let \[\mathcal{G}(\ell_1,\ell_2)=\{[0,\ell_1] \iso^+ [0,\ell_2]\}\] for $\ell_1,\ell_2>0$. 
\end{nota}

The sets $\mathcal{G}(\ell_1,\ell_2)$ are equipped with the $\Delta$-kelleyfication of the compact-open topology, which coincides with the compact-open topology, and with the pointwise convergence topology by \cite[Proposition~2.5]{Moore2}.

\bd \label{composition_map} Let $\gamma_1$ and $\gamma_2$ be two continuous maps from $[0,1]$ to some topological space such that $\gamma_1(1)=\gamma_2(0)$. The composite defined by 
\[
(\gamma_1 *_N \gamma_2)(t) = \begin{cases}
\gamma_1(2t)& \hbox{ if }0\leq t\leq \frac{1}{2},\\
\gamma_2(2t-1)& \hbox{ if }\frac{1}{2}\leq t\leq 1
\end{cases}
\]
is called the \textit{normalized composition}.
\ed

A \textit{multipointed $d$-space $X$} \cite{mdtop} is a triple $(|X|,X^0,\P^{top}X)$ where
\begin{itemize}[leftmargin=*]
	\item The pair $(|X|,X^0)$ is a multipointed space. The space $|X|$ is called the \textit{underlying space} of $X$ and the set $X^0$ the \textit{set of states} of $X$. 
	\item The set $\P^{top}X$ is a set of continuous maps from $[0,1]$ to $|X|$ called the \textit{execution paths}, satisfying the following axioms:
	\begin{itemize}
		\item For any execution path $\gamma$, one has $\gamma(0),\gamma(1)\in X^0$.
		\item Let $\gamma$ be an execution path of $X$. Then any composite $\gamma\phi$ with $\phi:[0,1]\iso^+ [0,1]$ is an execution path of $X$.
		\item Let $\gamma_1$ and $\gamma_2$ be two composable execution paths of $X$; then the normalized composition $\gamma_1 *_N \gamma_2$ is an execution path of $X$.
	\end{itemize}
\end{itemize}

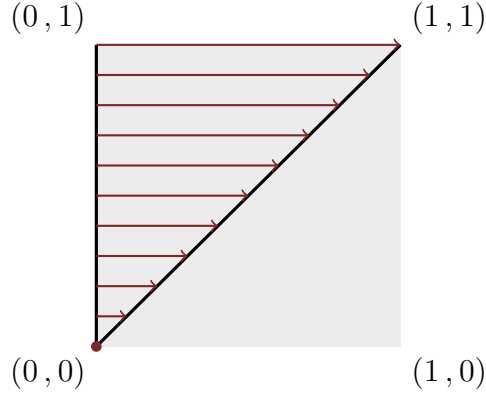
\begin{figure}
	\def\n{5}
	\begin{tikzpicture}[black,scale=4,pn/.style={circle,inner sep=0pt,minimum width=4pt,fill=dark-red}]
		\fill [color=gray!15] (0,0) -- (0,1) -- (1,1) -- (1,0) -- cycle;
		\draw[-] [very thick] (0,1) -- (0,0);
		\draw[-] [very thick] (0,0) -- (1,1);
		\foreach \n in {1,2,3,4,5,6,7,8,9,10}
		{\draw[dark-red][->][thick](0,\n/10) -- (\n/10,\n/10);}
		\draw (0,0) node[pn] {} node[black,below left] {$(0\,,0)$};
		\draw (1,1) node[black,above right] {$(1\,,1)$};
		\draw (0,1) node[black,above left] {$(0\,,1)$};
		\draw (1,0) node[black,below right] {$(1\,,0)$};
	\end{tikzpicture}
	\caption{$|X|=[0,1]\p [0,1]$, $X^0=\{0\}\p [0,1] \cup \{(x,x)\mid x\in [0,1]\}$, $\P^{top}_{(0,t),(t,t)}X= \mathcal{G}(1,1)$ for all $t\in ]0,1]$, $\P^{top}_{(0,0),(0,0)}X= \{(0,0)\}$ and $\P^{top}_{\alpha,\beta}X=\varnothing$ otherwise, there is no composable execution paths.}
	\label{contracting}
\end{figure}

\begin{rem}
	The subset $X^0$ is not necessarily discrete: for any topological space $Z$, the triple $(Z,Z,\varnothing)$ defines a multipointed $d$-space. Figure~\ref{contracting} depicts another example already given in \cite[Figure~1]{Moore2}.
\end{rem} 

A map $f:X\to Y$ of multipointed $d$-spaces is a map of multipointed spaces from $(|X|,X^0)$ to $(|Y|,Y^0)$ such that for any execution path $\gamma$ of $X$, the composite map $f. \gamma$ is an execution path of $Y$. 

\begin{nota}
	The category of multipointed $d$-spaces is denoted by $\ptop{\mathcal{G}}$. 
\end{nota}

The category of multipointed $d$-spaces is locally presentable by \cite[Theorem~3.5]{mdtop}.  Let us recall the three model structures of multipointed $d$-spaces.

\bth \label{three} \cite[Theorem~6.14]{QHMmodel} 
Let $r\in \{q,m,h\}$. There exists a unique model structure on $\ptop{\mathcal{G}}$ such that: 
\begin{itemize}[leftmargin=*]
	\item A map of multipointed $d$-spaces $f:X\to Y$ is a weak equivalence if and only if $f^0:X^0\to Y^0$ is a bijection and for all $(\alpha,\beta)\in X^0\p X^0$, the continuous map $\P^{top}_{\alpha,\beta}X\to \P^{top}_{f(\alpha),f(\beta)}Y$ is a weak equivalence of the r-model structure of $\top$.
	\item A map of flows $f:X\to Y$ is a fibration if and only if for all $(\alpha,\beta)\in X^0\p X^0$, the continuous map $\P^{top}_{\alpha,\beta}X\to \P^{top}_{f(\alpha),f(\beta)}Y$ is a fibration of the r-model structure of $\top$.
\end{itemize}
This model structure is accessible and all objects are fibrant. The q-model structure is even combinatorial.
\eth

\begin{nota}
	These three model categories are denoted by $\ptop{\mathcal{G}}_q$, $\ptop{\mathcal{G}}_m$ and $\ptop{\mathcal{G}}_h$.
\end{nota}

\begin{nota}
	The subset of execution paths from $\alpha$ to $\beta$ is the set of $\gamma\in\P^{top} X$  such that $\gamma(0)=\alpha$ and $\gamma(1)=\beta$; it is denoted by $\P^{top}_{\alpha,\beta} X$: $\alpha$ is called the \textit{initial state} and $\beta$ the \textit{final state} of such a $\gamma$. 
\end{nota}

The set $\P^{top}_{\alpha,\beta} X$ is equipped with the $\Delta$-kelleyfication of the relative topology induced by the inclusion $\P^{top}_{\alpha,\beta} X\subset \ttop([0,1],|X|)$.

\begin{nota}
	Let $[0] = \{()\}$ and $[n] = \{0,1\}^n$ for $n \geq 1$. By convention, one has $\{0,1\}^0=[0]=\{()\}$. The set $[n]$ is equipped with the product ordering $\{0<1\}^n$. Let $0_n=(0,\dots,0) \in \{0,1\}^n$ and $1_n=(1,\dots,1) \in \{0,1\}^n$
\end{nota}

Let $\delta_i^\alpha : [n-1] \rightarrow [n]$ be the \textit{coface map} defined for $1\leq i\leq n$ and $\alpha \in \{0,1\}$ by \[\delta_i^\alpha(\epsilon_1, \dots, \epsilon_{n-1}) = (\epsilon_1,\dots, \epsilon_{i-1}, \alpha, \epsilon_i, \dots, \epsilon_{n-1}).\] The small category $\square$, called sometimes the \textit{box category},  is by definition the subcategory of the category of sets with the set of objects $\{[n],n\geq 0\}$ and generated by the coface maps $\delta_i^\alpha$. They satisfy the cocubical relations $\delta_j^\beta \delta_i^\alpha = \delta_i^\alpha \delta_{j-1}^\beta $ for $i<j$ and for all $(\alpha,\beta)\in \{0,1\}^2$. If $p>q\geq 0$, then the set of morphisms $\square([p],[q])$ is empty. If $p = q$, then the set $\square([p],[p])$ is the singleton $\{\id_{[p]}\}$. For $0\leq p \leq q$, all maps of $\square$ from $[p]$ to $[q]$ are one-to-one. A good reference for presheaves is \cite{MR1300636}. 

\bd \cite{Brown_cube} The category of presheaves over $\square$, denoted by $\square^{op}\set$, is called the category of \textit{precubical sets}. 
\ed 

Let us expand the above definition. A precubical set $K$ consists of a family of sets $(K_n)_{n \geq 0}$ and of set maps $\de_i^\alpha:K_n \rightarrow K_{n-1}$ with $1\leq i \leq n$ and $\alpha\in\{0,1\}$ satisfying the cubical relations $\de_i^\alpha\de_j^\beta = \de_{j-1}^\beta \de_i^\alpha$ for any $\alpha,\beta\in \{0,1\}$ and for $i<j$. An element of $K_n$ is called a \textit{$n$-cube}. An element of $K_0$ is also called a vertex of $K$. A precubical set $K$ is of dimension $n\geq 0$ if $K_n\neq \varnothing$ and $K_p=\varnothing$ for $p>n$. 

\bd
	For all $n\geq 0$, the \textit{$n$-cube} is the precubical set $\square[n]=\square(-,[n])$. It is of dimension $n$.
\ed

There exists a functor $\square(K):(\square\ddownarrow K) \to \square^{op}\set$ (in the notation $(\square\ddownarrow K)$, every object $[n]$ of $\square$ is identified with the precubical set $\square[n]$) which takes the map of precubical sets $\square[n]\to K$ to $\square[n]$. It is a general property of presheaves that $K = \liminj \square(K)$, and the latter colimit is denoted by \[\liminj_{\square[n]\rightarrow K} \square[n].\] Let $\dim(x)=n$ if $x\in K_n$. Let \[K_{\leq n} = \liminj_{\substack{\square[p]\to K\\p\leq n}} \square[p].\] The \textit{boundary} of $\square[n]$ is the precubical set \[\de \square[n] = \square[n]_{\leq n-1}.\] In particular, one has $\de \square[0] = \varnothing$. The precubical set $\de \square[n]$ is of dimension $n-1$ for all $n\geq 1$.

\bd
A \textit{cocubical object} of a category $\C$ is a functor $\square\to \C$.
\ed

\begin{nota}
	Let $\C$ be a cocomplete category. Let $X:\square\to \C$ be a cocubical object of $\C$. It gives rise to a small diagram $(\square\ddownarrow K) \to \C$ for all precubical sets $K$. Its colimit is denoted by \[\widehat{X}(K)=\liminj_{\square[n]\to K} X([n]).\]
\end{nota}

\bp \cite[Proposition~2.3.2]{realization} \label{eq-cat}
Let $\C$ be a cocomplete category. The mapping $X\mapsto \widehat{X}$ induces an equivalence of categories between the category of cocubical objects of $\C$ and the colimit-preserving preserving functors from $\square^{op}\set$ to $\C$.
\ep

\bd \label{geom}
Let $K$ be a precubical set. The cocubical topological space defined on objects by $[n]\mapsto [0,1]^n$ and on maps by 
$(\delta_i^\alpha : [n-1] \rightarrow [n])\mapsto ((\epsilon_1, \dots, \epsilon_{n-1})\mapsto (\epsilon_1,\dots, \epsilon_{i-1}, \alpha, \epsilon_i, \dots, \epsilon_{n-1}))$ gives rise to a colimit-preserving functor from precubical sets to topological spaces denoted by 
\[
|K|_{geom} = \liminj_{\square[n]\to K} [0,1]^n.
\]
The space $|K|_{geom}$ is called the \textit{geometric realization} of $K$.
\ed

The topological space $|K|_{geom}$ is a CW-complex. It is equipped with the final topology which is $\Delta$-generated and Hausdorff.

\bd
Let $U$ be a topological space. A \textit{(Moore) path} in $U$ consists of a continuous map $\gamma:[0,\ell]\to U$ with $\ell>0$. The real number $\ell>0$ is called the \textit{length} of the path. It can be extended to a continuous map $\gamma:[0,+\infty[\to U$ such that $\gamma(t)=\gamma(\ell)$ for all $t\geq \ell$. The path $\gamma$ is a path from $\gamma(0)$ to $\gamma(\ell)$.
\ed

The notions of \textit{stop interval} and \textit{regular path} appear in \cite[Definition~1.1]{reparam}.

\bd \label{regular}
Let $U$ be a Hausdorff topological space. Let $\gamma:[0,\ell]\to U$ be a Moore path in $U$ with $\ell>0$. A \textit{stop interval} of $\gamma$ is an interval $[a,b]\subset [0,\ell]$ with $a<b$ such that the restriction $\gamma\rest_{[a,b]}$ is constant and such that $[a,b]$ is maximal for this property. The path $\gamma$ is \textit{regular} if it has no stop intervals.
\ed

Since $U$ is in this paper the geometric realization of a precubical set (cf. Definition~\ref{geom}), the Hausdorff condition in Definition~\ref{regular} is not restrictive. The Hausdorff condition is used in \cite{reparam} for two reasons: $U$ must have closed points and unique sequential limits. It is worth mentioning that the proof of \cite[Proposition~3.7]{reparam} is incorrect and that a correction is published in \cite{reparam-fixed}. Note also that the proof of \cite[Corollary~3.5]{reparam} is false, which has no consequence nowhere. Only two theorems from \cite{reparam} are used in this paper: 
\begin{itemize}[leftmargin=*]
	\item \cite[Proposition~2.22]{reparam}, which is a topological lemma involving only the segment $[0,1]$, in the proof of Proposition~\ref{viM}.
	\item \cite[Proposition~3.8]{reparam} in the remark before  Proposition~\ref{viM}. The proof of \cite[Proposition~3.8]{reparam} requires that the limit of a convergent sequence is unique.
\end{itemize}

\bd
Let $\gamma_1:[0,\ell_1]\to U$ and $\gamma_2:[0,\ell_2]\to U$ be two paths in the topological space $U$ such that $\gamma_1(\ell_1)=\gamma_2(0)$. The \textit{Moore composition} $\gamma_1*\gamma_2:[0,\ell_1+\ell_2]\to U$ is the Moore path defined by 
\[
(\gamma_1*\gamma_2)(t)=
\begin{cases}
\gamma_1(t) & \hbox{ for } t\in [0,\ell_1]\\
\gamma_2(t-\ell_1) &\hbox{ for }t\in [\ell_1,\ell_1+\ell_2].
\end{cases}
\]
The Moore composition of Moore paths is strictly associative. The Moore composition of two regular paths of a Hausdorff space is regular.
\ed

\bd \label{def-dpath}
Let $n\geq 0$. A \textit{(Moore) $d$-path} of length $\ell\geq 0$ of $\square[n]$ (or in the geometric realization of $\square[n]$) is a continuous map $\gamma:[0,\ell]\to [0,1]^n$ with $\ell\geq 0$ which is nondecreasing with respect to each axis of coordinates. It is \textit{tame} if $\gamma(0),\gamma(\ell)\in \{0,1\}^n$. Let $K$ be a precubical set. Let $c\in K_n$ with $n\geq 0$ be an $n$-cube of $K$. A \textit{$d$-path} of $c$ of length $\ell \geq 0$ is a composite continuous map denoted by $[c;\gamma]:[0,\ell] \to |K|_{geom}$ such that $\gamma:[0,\ell]\to [0,1]^n$ is a $d$-path of length $\ell\geq 0$ with $[c;\gamma]=|c|_{geom}\gamma$. A \textit{$d$-path} of $K$ (or in the geometric realization of $K$) of length $\ell\geq 0$ is a continuous path $[0,\ell] \to |K|_{geom}$ which is a Moore composition of the form $[c_1;\gamma_1] * \dots *[c_p;\gamma_p]$ with $p\geq 1$. The $d$-path $\gamma$ is \textit{tame} if each $\gamma_i$ for $1\leq i\leq p$ is tame. The point $\gamma(0)$ is called the \textit{initial state} of $\gamma$ and the point $\gamma(\ell)$ is called the \textit{final state} of $\gamma$. We refer e.g. to \cite[Section~2.1 and Section~2.5]{MR4070250} for further details. 
\ed 

\begin{figure}
	\begin{tikzpicture}
		\draw (0,0) -- (0,2) -- (2,2) -- (2,0) -- (0,0);
		\draw (2,2) -- (4,2) -- (4,0) -- (2,0);
		\draw (2,2) -- (2,4) -- (4,4) -- (4,2);
		\draw (5,0) -- (5,2) -- (7,2) -- (7,0) -- (5,0);
		\draw (7,2) -- (9,2) -- (9,0) -- (7,0);
		\draw (7,2) -- (7,4) -- (9,4) -- (9,2);	
		\draw[->, line width=0.4mm, color=dark-red, smooth] plot coordinates {(0,0) (0.5,1.5) (2,2)};
        \draw[->, line width=0.4mm, color=dark-red, smooth] plot coordinates {(2,2) (3.5,2.5) (4,4)};		
		\draw[->, line width=0.4mm, color=dark-red, smooth] plot coordinates {(5,0) (6,0.7) (7,1)};
		\draw[->, line width=0.4mm, color=dark-red, smooth] plot coordinates {(7,1) (8,1.3) (8.2,2)};
		\draw[->, line width=0.4mm, color=dark-red, smooth] plot coordinates {(8.2,2) (8.5,3) (9,4)};	
	\end{tikzpicture}
	\caption{The left $d$-path is tame, the right $d$-path is not tame.}
	\label{tame-nontame}
\end{figure}
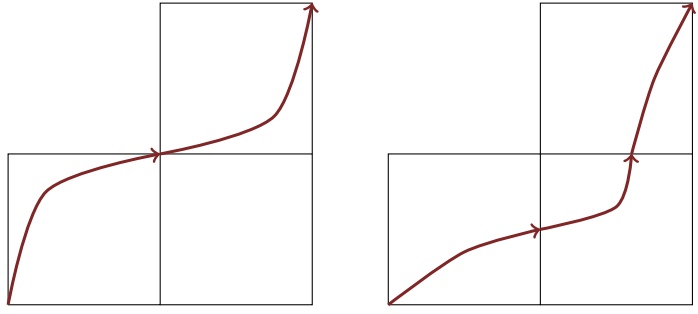

Figure~\ref{tame-nontame} depicts an example of a tame $d$-path and of a non-tame $d$-path: the exit points of the right $d$-path from the bottom left square and from the bottom right square are not vertices indeed. Figure~\ref{nontame} gives another example of a non-tame $d$-path in the boundary of the $3$-cube. Note that all $d$-paths in the boundary of the $2$-cube are tame.

\begin{nota}
	The \textit{set} of \textit{all} $d$-paths of length $1$ of a precubical set $K$ from $\alpha\in K_0$ to $\beta\in K_0$ is denoted by $\vec{P}_{\alpha,\beta}(K)$.
\end{nota}

Now we can introduce the \textit{tame regular realization} of a precubical set $K$ as follows. The \textit{regular realization} of a precubical set is defined in Section~\ref{question-r}.

\bp \label{cocubical-reg-cube}
Let $n\geq 1$. The following data assemble into a multipointed $d$-space $|\square[n]|^t_{reg}$ called the \textit{regular $n$-cube}:
\begin{itemize}
	\item The underlying space is the topological $n$-cube $[0,1]^n$.
	\item The set of states is $\{0,1\}^n\subset [0,1]^n$.
	\item The set of execution paths from $\underline{a}$ to $\underline{b}$ with $\underline{a}<\underline{b} \in \{0<1\}^n$ is the set of regular $d$-paths $[0,1]\to [0,1]^n$ from $\underline{a}$ to $\underline{b}$.
	\item The set of execution paths from $\underline{a}$ to $\underline{b}$ with $\underline{a}\geq \underline{b}$ is empty.
\end{itemize}
Let $|\square[0]|^t_{reg} = \{()\}$. 	Let $\delta_i^\alpha : [0,1]^{n-1} \rightarrow [0,1]^n$ be the continuous map defined for $1\leq i\leq n$ and $\alpha \in \{0,1\}$ by $\delta_i^\alpha(\epsilon_1, \dots, \epsilon_{n-1}) = (\epsilon_1,\dots, \epsilon_{i-1}, \alpha, \epsilon_i, \dots, \epsilon_{n-1})$. By convention, let $[0,1]^0=\{()\}$. The mapping $[n]\mapsto |\square[n]|^t_{reg}$ yields a well-defined cocubical multipointed $d$-space.
\ep

\bpf
The multipointed $d$-space $|\square[n]|^t_{reg}$ is well defined for all $n\geq 0$ because the normalized composition of two regular paths is regular and because the reparametrization of a regular path by a homeomorphism is regular.
\epf

Using Proposition~\ref{eq-cat}, this leads to the definition: 

\bd \label{treg_rea}
Let $K$ be a precubical set. The \textit{tame regular realization} of $K$ is the multipointed $d$-space 
\[
|K|^t_{reg} = \liminj_{\square[n]\rightarrow K} |\square[n]|^t_{reg}.
\]
This yields a colimit-preserving functor from precubical sets to multipointed $d$-spaces.
\ed

Since the underlying space functor from multipointed $d$-spaces to spaces is colimit-preserving, the topological space $|K|_{geom}$ is the underlying space of the multipointed $d$-space $|K|^t_{reg}$. By \cite[Proposition~6.5]{QHMmodel}, the forgetful functor $X\mapsto (|X|,X^0)$ from multipointed $d$-spaces to multipointed spaces which forgets the set of execution paths is topological in the sense of \cite[Definition~21.1]{topologicalcat}. This means that colimits are calculated by a final structure. Consequently, the execution paths of $|K|^t_{reg}$ are exactly the nonconstant tame  regular $d$-paths of $K$.

\section{Homotopy theory of Moore flows}
\label{qhm_sec}

This section starts by the necessary reminders about flows and Moore flows for the readers not familiar with this topic. The end of the section introduces the q-model structure, the h-model structure and the m-model structure of Moore flows, as generalizations of results expounded in \cite{QHMmodel} for flows: the proofs are very similar indeed.

\bd \cite[Definition~4.11]{model3} \label{def-flow}
A \textit{flow} is a small semicategory enriched over the closed monoidal category $(\top,\p)$. The corresponding category is denoted by $\dtop$. 
\ed

Let us expand the above definition. A \textit{flow} $X$ consists of a topological space $\P X$ of \textit{execution paths}, a discrete space $X^0$ of \textit{states}, two continuous maps $s$ and $t$ from $\P X$ to $X^0$ called the source and target map respectively, and a continuous and associative map $*:\{(x,y)\in \P X\p \P X; t(x)=s(y)\}\longrightarrow \P X$ such that $s(x*y)=s(x)$ and $t(x*y)=t(y)$. Let $\P_{\alpha,\beta}X = \{x\in \P X\mid s(x)=\alpha \hbox{ and } t(x)=\beta\}$: it is the space of execution paths from $\alpha$ to $\beta$, $\alpha$ is called the initial state and $\beta$ is called the final state. Note that the composition is denoted by $x*y$, not by $y\circ x$. The category $\dtop$ is locally presentable by \cite[Theorem~6.11]{Moore1}. 

\begin{nota} \label{basic-example-flow}
	A set can be viewed as a flow with an empty path space between each pair of states. The flow $\vI$ consists of the semicategory $0\rightarrow 1$ without identity morphisms.
\end{nota}

Let $X$ be a multipointed $d$-space. Consider for every $(\alpha,\beta)\in X^0 \p X^0$ the coequalizer of spaces \[\P_{\alpha,\beta}X = \liminj\left( \P^{top}_{\alpha,\beta}X\p \mathcal{G}(1,1) \rightrightarrows \P^{top}_{\alpha,\beta}X\right)\] where the two maps are $(c,\phi)\mapsto c$ and $(c,\phi)\mapsto c.\phi$. Let $[-]_{\alpha,\beta}:\P^{top}_{\alpha,\beta}X \rightarrow \P_{\alpha,\beta}X$ be the canonical map.

\bth \cite[Theorem~7.2]{mdtop} \label{decomposing}
Let $X$ be a multipointed $d$-space. Then there exists a flow $\Cat(X)$ with $\Cat(X)^0=X^0$, $\P_{\alpha,\beta}\Cat(X)= \P_{\alpha,\beta}X$ and the composition law $*:\P_{\alpha,\beta}X \p \P_{\beta,\gamma}X \rightarrow \P_{\alpha,\gamma}X$ is for every triple $(\alpha,\beta,\gamma)\in X^0\p X^0\p X^0$ the unique map making the following diagram commutative:
\[
\xymatrix@C=4em@R=4em{
	\P^{top}_{\alpha,\beta}X \p \P_{\beta,\gamma}^{\mathcal{G}}X
	\fr{*_{N}}\fd{[-]_{\alpha,\beta}\p [-]_{\beta,\gamma}} &
	\P_{\alpha,\gamma}^{\mathcal{G}}X \fd{[-]_{\alpha,\gamma}} \\
	\P_{\alpha,\beta}X \p \P_{\beta,\gamma}X \fr{} &
	\P_{\alpha,\gamma}X}
\] 
where $*_N$ is the normalized composition (cf. Definition~\ref{composition_map}).  The mapping $X \mapsto cat(X)$ induces a functor from $\ptop{\mathcal{G}}$ to $\dtop$.  \eth

\bd \label{cat-func} The functor $\Cat:\ptop{\mathcal{G}}\to \dtop$ is called the \textit{categorization functor}. \ed

\begin{nota} \label{defGcat}
	The enriched small category $\mathcal{G}$ is defined as follows: 
\begin{itemize}[leftmargin=*]
	\item The set of objects is the open interval $]0,+\infty[$.
	\item The space of maps from $\ell_1$ to $\ell_2$ is the space $\mathcal{G}(\ell_1,\ell_2)$ defined in Notation~\ref{defG}.
	\item For every $\ell_1,\ell_2,\ell_3>0$, the composition map \[\mathcal{G}(\ell_1,\ell_2)\p \mathcal{G}(\ell_2,\ell_3) \to \mathcal{G}(\ell_1,\ell_3)\] is induced by the composition of continuous maps. 
\end{itemize}
\end{nota}

The enriched category $\mathcal{G}$ is an example of a \textit{reparametrization category} in the sense of Definition~\ref{def-reparam} which is different from the terminal category. It is introduced in \cite[Proposition~4.9]{Moore1}. 

By \cite[Proposition~3.8]{reparam}, the topological space $\P_{\alpha,\beta}X$ of Theorem~\ref{decomposing} is exactly the space of traces in the sense of \cite{MR2521708} when $|X|$ is Hausdorff and all execution paths of the multipointed $d$-space $X$ from $\alpha$ to $\beta$ are regular. This happens e.g. when $X$ is a cellular object of the q-model category of multipointed $d$-spaces of Theorem~\ref{three}. It is not true for general multipointed $d$-spaces as shown by Proposition~\ref{viM}. 

To have such a fact for all multipointed $d$-spaces, the enriched category $\mathcal{G}$ must be replaced by the enriched category $\mathcal{M}$ of \cite[Proposition~4.11]{Moore1} in the definition of a multipointed $d$-space. The enriched category $\mathcal{M}$ is another example of reparametrization category in the sense of Definition~\ref{def-reparam}. It is out of the scope of this paper which treats only the case of \textit{regular} $d$-paths. The reader might be interested in \cite{Moore3} to read the first results in this direction.

\bp \label{viM}
Consider the multipointed $d$-space $\vI^{\mathcal{M}}$ such that
\begin{itemize}
	\item The underlying topological space of $\vI^{\mathcal{M}}$ is the segment $[0,1]$.
	\item The set of states of $\vI^{\mathcal{M}}$ is $\{0,1\}$.
	\item The set of execution paths of $\vI^{\mathcal{M}}$ from $0$ to $1$ is the set $\mathcal{M}(1,1)$ of nondecreasing surjective maps from $[0,1]$ to itself. There are no other execution paths in $\vI^{\mathcal{M}}$.
	\item There is no composition law.
\end{itemize}
Then $\P_{0,1}\vI^{\mathcal{M}}$ is not a singleton. 
\ep

\bpf
The space $\P_{0,1}\vI^{\mathcal{M}}$ is the quotient of the space $\mathcal{M}(1,1)$ by the action of the space $\mathcal{G}(1,1)$. By \cite[Proposition~2.5]{Moore2}, the compact-open topology of $\mathcal{G}(1,1)$ is $\Delta$-generated, and by \cite[Proposition~2.6]{Moore3}, the compact-open topology of $\mathcal{M}(1,1)$ is $\Delta$-generated as well (it is the same argument); therefore we can use the results of \cite{reparam} which is written with the compact-open topology. By \cite[Proposition~2.22]{reparam}, the quotient of $\mathcal{M}(1,1)$ by the action of $\mathcal{G}(1,1)$ is in bijection with the set of countable subsets of $[0,1]$. We deduce that $\P_{0,1}\vI^{\mathcal{M}}$ is not a singleton.
\epf

This implies that the flow $\Cat(\vI^{\mathcal{M}})$ is not equal to $\vI$ (see Notation~\ref{basic-example-flow}). Since the inclusion map $\mathcal{G}(1,1) \subset \mathcal{M}(1,1)$ is a homotopy equivalence, a q-cofibrant replacement of $\vI^{\mathcal{M}}$ is the multipointed $d$-space $\vI^{\mathcal{G}}$ such that
\begin{itemize}
	\item The underlying topological space of $\vI^{\mathcal{G}}$ is the segment $[0,1]$.
	\item The set of states of $\vI^{\mathcal{G}}$ is $\{0,1\}$.
	\item The set of execution paths of $\vI^{\mathcal{G}}$ from $0$ to $1$ is the set $\mathcal{G}(1,1)$. There are no other execution paths in $\vI^{\mathcal{G}}$.
	\item There is no composition law.
\end{itemize}
This phenomenon is general: the q-cofibrant replacement functor of $\ptop{\mathcal{G}}$ turns out to be a way of removing all non-regular execution paths of a multipointed $d$-space by preserving its causal structure.

\bd \label{gspace}
The enriched category of enriched presheaves from $\mathcal{G}$ to $\top$ is denoted by $\topdgr$. The underlying set-enriched category of enriched maps of enriched presheaves is denoted by $\topdgr_0$. The objects of $\topdgr_0$ are called the \textit{$\mathcal{G}$-spaces}. Let \[\mathbb{F}^{\mathcal{G}^{op}}_{\ell}U=\mathcal{G}(-,\ell)\p U \in \topdgr_0\] where $U$ is a topological space and where $\ell>0$.
\ed

The category $\topdgr_0$ is locally presentable by \cite[Proposition~5.1]{dgrtop}. 

\bp\label{ev-adj} \cite[Proposition~5.3 and Proposition~5.5]{dgrtop}
The category of enriched presheaves $\topdgr_0$ is a full reflective and coreflective subcategory of $\top^{\mathcal{G}^{op}_0}$.  For every $\mathcal{G}$-space $F:\mathcal{G}^{op}\to \top$, every $\ell>0$ and every topological space $X$, we have the natural bijection of sets \[\topdgr_0(\mathbb{F}^{\mathcal{G}^{op}}_{\ell}X,F) \iso \top(X,F(\ell)).\] 
\ep

\bth (\cite[Theorem~5.14]{Moore1}) \label{closedsemimonoidal}
Let $D$ and $E$ be two $\mathcal{G}$-spaces. Let 
\[
D \ot E = \int^{(\ell_1,\ell_2)} \mathcal{G}(-,\ell_1+\ell_2) \p D(\ell_1) \p E(\ell_2).
\]
The pair $(\topdgr_0,\ot)$ has the structure of a biclosed semimonoidal category.
\eth

\bp \label{Ftenseur} (\cite[Proposition~5.16]{Moore1})
Let $U,U'$ be two topological spaces. Let $\ell,\ell'>0$. There is the natural isomorphism of $\mathcal{G}$-spaces 
\[
\mathbb{F}^{\mathcal{G}^{op}}_{\ell}U \ot \mathbb{F}^{\mathcal{G}^{op}}_{\ell'}U' \iso \mathbb{F}^{\mathcal{G}^{op}}_{\ell+\ell'}(U\p U').
\]
\ep

\bp \label{tensor-product} \cite[Proposition~5.18]{Moore1}
Let $D$ and $E$ be two $\mathcal{G}$-spaces. Then there is a natural homeomorphism 
\[
\liminj (D \ot E) \iso \liminj D \p \liminj E.
\]
\ep

By \cite[Theorem~6.5(ii)]{MoserLyne}, since all topological spaces are fibrant and since $(\top,\p,\{0\})$ is a locally presentable base by \cite[Corollary~3.3]{dgrtop}, the category of $\mathcal{G}$-spaces $\topdgr_0$ can be endowed with the projective model structure associated with one of the three model structures $\top_q,\top_m,\top_h$. They are called the projective q-model structure (m-model structure, h-model structure resp.) and denoted by $\topdgrq$ ($\topdgrm$, $\topdgrh$ resp.). The three model structures are accessible. The fibrations are the objectwise fibrations of the corresponding model structure of $\top$. They are called the projective q-fibrations (projective m-fibrations, projective h-fibrations resp.). The weak equivalences are the objectwise weak equivalence of the corresponding model structure of $\top$. Since the projective m-fibrations of spaces are the projective h-fibrations of spaces, and since the weak equivalences of the projective m-model structure of $\mathcal{G}$-spaces are the weak equivalences of the projective q-model structure of $\mathcal{G}$-spaces, it implies that the projective m-model structure of $\mathcal{G}$-spaces is the mixing in the sense of \cite[Theorem~2.1]{mixed-cole} of the projective q-model structure of $\mathcal{G}$-spaces and of the projective h-model structure of $\mathcal{G}$-spaces. All $\mathcal{G}$-spaces are fibrant for these three model structures. 

When the reparametrization category $\mathcal{G}$ is replaced by the terminal category, these model structures coincide with the q-model structure $\top_q$, the m-model structure $\top_m$ and the h-model structure $\top_h$ respectively.

\bd \cite[Definition~6.2]{Moore1} \label{def-Moore-flow}
A \textit{Moore flow} is a small semicategory enriched over the biclosed semimonoidal category $(\topdgr_0,\ot)$ of Theorem~\ref{closedsemimonoidal}. The corresponding category is denoted by $\dtopG$. 
\ed

A Moore flow $X$ consists of a \emph{set of states} $X^0$, for each pair $(\alpha,\beta)$ of states a $\mathcal{G}$-space $\P_{\alpha,\beta}X$ of $\topdgr_0$ and for each triple $(\alpha,\beta,\gamma)$ of states an associative composition law \[*:\P_{\alpha,\beta}X \ot \P_{\beta,\gamma}X \to \P_{\alpha,\gamma}X.\] A map of Moore flows $f$ from $X$ to $Y$ consists of a set map \[f^0:X^0 \to Y^0\] (often denoted by $f$ as well if there is no possible confusion) together for each pair of states $(\alpha,\beta)$ of $X$ with a natural transformation \[\P f:\P_{\alpha,\beta}X \longrightarrow \P_{f(\alpha),f(\beta)}Y\] compatible with the composition. The topological space $\P_{\alpha,\beta}X(\ell)$ is denoted by $\P_{\alpha,\beta}^\ell X$ and is called the space of \emph{execution paths of length $\ell$}. A set can be viewed as a Moore flow with an empty $\mathcal{G}$-space of execution paths between each pair of states.

Note that by replacing $\mathcal{G}$ by the terminal category, we recover Definition~\ref{def-flow}.

\begin{nota}
	Let $D:\mathcal{G}^{op}\to \top$ be a $\mathcal{G}$-space. We denote by $\glob(D)$ the Moore flow defined as follows: 
	\[
	\begin{aligned}
		&\glob(D)^0 = \{0,1\}\\
		&\P_{0,0}\glob(D)=\P_{1,1}\glob(D)=\P_{1,0}\glob(D)=\varnothing\\
		&\P_{0,1}\glob(D)=D.
	\end{aligned}	
	\]
	There is no composition law. This construction yields a functor \[\glob:\topdgr_0\to \dtopG.\]
\end{nota}

The category $\dtopG$ is locally presentable by \cite[Theorem~6.11]{Moore1}. Let $X$ be a multipointed $d$-space. Let $\P^\ell_{\alpha,\beta}X$ be the subspace of continuous maps from $[0,\ell]$ to $|X|$ defined by $\P^\ell_{\alpha,\beta}X = \{t\mapsto \gamma\mu_\ell\mid \gamma\in \P^{top}_{\alpha,\beta}X\}$. By \cite[Theorem~4.12]{Moore2}, there exists a Moore flow $\moore^{\mathcal{G}}(X)$ such that: 
\begin{itemize}
	\item The set of states $X^0$ of $X$.
	\item For all $\alpha,\beta\in X^0$ and all real numbers $\ell>0$, $\P_{\alpha,\beta}^{\ell}\moore^{\mathcal{G}}(X) = \P_{\alpha,\beta}^{\ell}X$.
	\item For all maps $[0,\ell]\iso^+[0,\ell']$, a map $f:[0,\ell']\to |X|$ of $\P_{\alpha,\beta}^{\ell'}\moore^{\mathcal{G}}(X)$ is mapped to the map $[0,\ell]\iso^+[0,\ell']\stackrel{f}\to |X|$ of $\P_{\alpha,\beta}^{\ell}\moore^{\mathcal{G}}(X)$.
	\item For all $\alpha,\beta,\gamma\in X^0$ and all real numbers $\ell,\ell'>0$, the composition maps \[*:\P_{\alpha,\beta}^{\ell}\moore^{\mathcal{G}}(X) \p \P_{\beta,\gamma}^{\ell'}\moore^{\mathcal{G}}(X) \to \P_{\alpha,\gamma}^{\ell+\ell'}\moore^{\mathcal{G}}(X)\] is the Moore composition.
\end{itemize}

\bth \label{MG} (\cite[Theorem~4.12 and Appendix~B]{Moore2})
The mapping above induces a functor $\moore^{\mathcal{G}}:\ptop{\mathcal{G}}\to\dtopG$ which is a right adjoint.
\eth

\begin{nota}
	Denote by $\lmoore^{\mathcal{G}}:\dtopG\to\ptop{\mathcal{G}}$ the left adjoint of $\moore^{\mathcal{G}}:\ptop{\mathcal{G}}\to \dtopG$. 
\end{nota}

Consider a Moore flow $X$. For all $\alpha,\beta\in X^0$, let $X_{\alpha,\beta}=\liminj \P_{\alpha,\beta}X$. Let $(\alpha,\beta,\gamma)$ be a triple of states of $X$. The composition law of the Moore flow $X$ induces, using Proposition~\ref{tensor-product}, a continuous map \[X_{\alpha,\beta} \p X_{\beta,\gamma} \iso \liminj(\P_{\alpha,\beta}X \ot \P_{\beta,\gamma}X) \longrightarrow \liminj \P_{\alpha,\gamma}X \iso X_{\alpha,\gamma}\] which is associative in an obvious sense. We obtain a flow $\lmoore(X)$ such that:
\begin{itemize}
	\item  The set of states is $X^0$.
	\item For all $\alpha,\beta\in X^0$, one has  $\P_{\alpha,\beta}\lmoore X=X_{\alpha,\beta}$.
	\item For all $\alpha,\beta,\gamma\in X^0$, the composition law is the map $X_{\alpha,\beta}\p X_{\beta,\gamma}\to X_{\alpha,\gamma}$ above defined.
\end{itemize}
This construction yields a well-defined functor \[\lmoore:\dtopG \longrightarrow \dtop\] by \cite[Proposition~10.5]{Moore1} which is a left adjoint by \cite[Theorem~10.7]{Moore1}. 

\begin{nota}
	Denote by $\moore:\dtop\to \dtopG$ a right adjoint of the functor $\lmoore:\dtopG \longrightarrow \dtop$.
\end{nota}

We now recall the theorem:

\bth \label{decomposing2} (\cite[Theorem~8.11]{Moore2}) There is the isomorphism of functors \[\Cat\iso\lmoore\moore^{\mathcal{G}}\] 
where $\Cat:\ptop{\mathcal{G}}\to \dtop$ from multipointed $d$-spaces to flows is the functor of Definition~\ref{cat-func}. \eth

We want to conclude this section by constructing the \{q,h,m\}-model structures of Moore flows.

\bd A \textit{$\topdgr_0$-graph} $X$ consists of a pair \[(X^0,(\P_{\alpha,\beta}X)_{(\alpha,\beta)\in X^0\p X^0})\] such that $X^0$ is a set and such that each $\P_{\alpha,\beta}X$ is a $\mathcal{G}$-space. A map of $\topdgr_0$-graphs $f:X\to Y$ consists of a set map $f^0:X^0\to Y^0$ (called the \textit{underlying set map}) together with a map $\P_{\alpha,\beta}X\to \P_{f^0(\alpha),f^0(\beta)}Y$ of $\mathcal{G}$-spaces for all $(\alpha,\beta)\in X^0\p X^0$. The composition is defined in an obvious way. The corresponding category is denoted by $\pf(\topdgr_0)$. 
\ed

\bth \label{qhmMooreFlow}
Let $(\C,\F,\W)$ be either the projective q-model structure, or the projective m-model structure, or the projective h-model structure of $\mathcal{G}$-spaces. Then the category of Moore flows can be endowed with an accessible model structure characterized as follows: 
\begin{itemize}[leftmargin=*]
	\item A map of Moore flows $f:X\to Y$ is a weak equivalence if and only if $f^0:X^0\to Y^0$ is a bijection and $\P f:\P_{\alpha,\beta}X\to \P_{f(\alpha),f(\beta)}Y$ belongs to $\W$ for all $(\alpha,\beta)\in X^0\p X^0$.
	\item A map of Moore flows $f:X\to Y$ is a fibration if and only if $\P f:\P_{\alpha,\beta}X\to \P_{f(\alpha),f(\beta)}Y$ belongs to $\F$ for all $(\alpha,\beta)\in X^0\p X^0$.
\end{itemize}
All objects are fibrant. These three model structures are called the \textit{q-model structure}, the \textit{m-model structure} and the \textit{h-model structure} of $\dtopG$. The m-model structure of Moore flows is the mixing in the sense of \cite[Theorem~2.1]{mixed-cole} of the q-model structure and the h-model structure of Moore flows. The q-model structure coincides with the one of \cite[Theorem~8.8]{Moore1}. Every q-cofibration of Moore flows is an m-cofibration of Moore flows and every m-cofibration of Moore flows is an h-cofibration of Moore flows.  
\eth

\begin{rem}
	Theorem~\ref{qhmMooreFlow} is already proved in \cite[Theorem~7.4]{QHMmodel} if the reparametrization category $\mathcal{G}$ is replaced by the terminal category, which enables us to obtain the \textit{q-model structure}, the \textit{m-model structure} and the \textit{h-model structure} of $\dtop$. All these model structures are accessible. The q-model structure is even combinatorial.
\end{rem}

\bpf[Sketch of proof]
Choose $(\C,\F,\W)$. By \cite[Corollary~5.6]{QHMmodel}, there exists a unique model structure on the category $\pf(\topdgr_0)$ such that the weak equivalences are the maps of enriched graphs which induce a bijection between the sets of vertices and which are objectwise weak equivalences and such that the fibrations are the objectwise fibrations. This model structure is accessible and all enriched graphs are fibrant. There is a forgetful functor \[\pf:\dtopG\longrightarrow\pf(\topdgr_0)\] which forgets the composition law. This functor is a right-adjoint. The left adjoint is the free Moore flow generated by the enriched graph. The semicategorical description of Moore flows is crucial: see \cite[Proposition~7.3]{QHMmodel} for a description of the left adjoint. Like in \cite[Theorem~7.4]{QHMmodel}, the model structure on enriched graphs of $\mathcal{G}$-spaces can be right-lifted along this right adjoint by using \cite[Theorem~2.1]{QHMmodel}, namely the Quillen path object argument in model categories where all objects are fibrant, with the path functor $\cocyl:\dtopG\to\dtopG$ defined on objects by $\cocyl(X)^0=X^0$, and for all $(\alpha,\beta)\in X^0\p X^0$ and for all $\ell>0$ by $\P_{\alpha,\beta}^\ell\cocyl(X) = \ttop([0,1],\P_{\alpha,\beta}^\ell X)$ with an obvious definition of the composition law. The hypotheses of \cite[Theorem~2.1]{QHMmodel} are satisfied because they are satisfied objectwise. Indeed, for all $(\alpha,\beta)\in X^0\p X^0$ and for all $\ell>0$, the diagonal of the topological space $\P_{\alpha,\beta}^\ell X$ factors as a composite \[\P_{\alpha,\beta}^\ell X \longrightarrow \ttop([0,1],\P_{\alpha,\beta}^\ell X) \longrightarrow \P_{\alpha,\beta}^\ell X\p \P_{\alpha,\beta}^\ell X\] where the left-hand map is a homotopy equivalence and the right-hand map is a r-fibration for $r\in\{q,m,h\}$. The q-model structure of this theorem coincides with the one of \cite[Theorem~8.8]{Moore1} because the two model structures have the same class of fibrations and weak equivalences. The last sentence is a general fact about mixed model structures (see \cite[Theorem~2.1 and Proposition~3.6]{mixed-cole}).
\epf

\begin{nota}
	These three model categories are denoted by $\dtopG_q$, $\dtopG_m$ and $\dtopG_h$.
\end{nota}

\begin{rem}
	The maps $C:\varnothing \to \{0\}$ and $R:\{0,1\}\to \{0\}$ are cofibrations of Moore flows. See \cite{Nonunital} for some explanations about the importance of the cofibration $R:\{0,1\}\to \{0\}$.
\end{rem}

By equipping the category of multipointed $d$-spaces with its q-model structure (cf. Theorem~\ref{three}), the functor $\moore^{\mathcal{G}}:\ptop{\mathcal{G}}\to\dtopG$ becomes a right Quillen equivalence between the q-model structures by \cite[Theorem~8.1]{Moore2}. Moreover, the unit and counit of the Quillen adjunction $\lmoore^{\mathcal{G}} \dashv \moore^{\mathcal{G}}$ induce isomorphisms for the q-cofibrant objects \cite[Theorem~7.6 and Corollary~7.9]{Moore2}. The functor $\lmoore:\dtopG \longrightarrow \dtop$ yields a left Quillen equivalence between the q-model structures by \cite[Theorem~10.9]{Moore1}. Finally, by \cite[Theorem~8.8]{Moore2}, the functor $\Cat\iso\lmoore\moore^{\mathcal{G}}:\ptop{\mathcal{G}}\to \dtop$ is neither a left adjoint nor a right adjoint. However, by \cite[Theorem~8.14]{Moore2}, it takes q-cofibrant multipointed $d$-spaces to q-cofibrant flows and its total left derived functor in the sense of \cite{HomotopicalCategory}, namely $X\mapsto \Cat(X^{cof})$ where $(-)^{cof}$ is a q-cofibrant replacement of the q-model structure of multipointed $d$-spaces, induces an equivalence of categories between the homotopy categories of the q-model structures. An inverse functor of $\Cat:\ptop{\mathcal{G}}\to \dtop$ up to homotopy is the functor $X\mapsto \lmoore^{\mathcal{G}}(\moore(X)^{cof})$ where $(-)^{cof}$ is now a q-cofibrant replacement of the q-model structure of Moore flows. By Proposition~\ref{viM}, one has $\Cat(\vI^{\mathcal{M}})\neq \vI$. This is not contradictory. Indeed, the q-cofibrant replacement of the multipointed $d$-space $\vI^{\mathcal{M}}$ is the multipointed $d$-space $\vI^{\mathcal{G}}$ described in the comment following Proposition~\ref{viM} and $\Cat(\vI^{\mathcal{G}})= \vI$.

It is worth noting that the h-model structures of Moore flows and of flows, as well as the h-model structure of multipointed $d$-spaces recalled in Theorem~\ref{three} do not coincide with the Hurewicz model structure given by \cite[Corollary~5.23]{Barthel-Riel}. This one exists as well for Moore flows, flows, and multipointed $d$-spaces because all these categories satisfy the monomorphism hypothesis of \cite[Definition~5.16]{Barthel-Riel}, being locally presentable, and because all of them are enriched, tensored and cotensored over ($\Delta$-Hausdorff or not) $\Delta$-generated spaces. The proof of the latter fact is similar to the proof that they are simplicial, tensored and cotensored. The proof for flows is written in \cite[Section~3.3]{realization}. The proof for multipointed $d$-spaces is sketched in \cite[Appendix~B]{mdtop}. The proof for Moore flows is left to the reader. The reason of this similarity is that every $\Delta$-generated space is homeomorphic to the disjoint sum of its path-connected components by \cite[Proposition~2.8]{mdtop}, exactly like simplicial sets. These ``genuine'' Hurewicz model structures provided by \cite[Corollary~5.23]{Barthel-Riel} are not used in this paper. In fact, by now, there are no known applications of them in the theory of flows, Moore flows or multipointed $d$-spaces.

\section{\mins{L_1}-arc length of \mins{d}-paths of precubical sets}
\label{length_sec}

\bd \label{path-with-length}
Let $K$ be a precubical set. Let $(\alpha,\beta)\in K_0\p K_0$. Let $\ell>0$. Let $\vec{R}^\ell_{\alpha,\beta}(K)$ be the subspace of continuous maps from $[0,\ell]$ to $|K|_{geom}$ defined by \[\vec{R}^\ell_{\alpha,\beta}(K) = \{t\mapsto \gamma\mu_\ell\mid \gamma\in \P^{top}_{\alpha,\beta}|K|^t_{reg}\}.\] 
Its elements are called \textit{the (nonconstant) tame regular $d$-paths of $K$ (of length $\ell$)} from $\alpha$ to $\beta$. Let 
\[
\vec{R}^\ell(K) = \coprod_{(\alpha,\beta)\in K_0\p K_0} \vec{R}^\ell_{\alpha,\beta}(K).
\]
\ed

The definition of $\vec{R}^\ell_{\alpha,\beta}(K)$ is not restrictive. Indeed, we have: 

\bp Let $K$ be a precubical set. Let $\phi:[0,\ell]\iso^+ [0,\ell]$. Let $\gamma\in \vec{R}^\ell_{\alpha,\beta}(K)$. Then $\gamma\phi\in \vec{R}^\ell_{\alpha,\beta}(K)$. \ep

\bpf By definition of $\vec{R}^\ell_{\alpha,\beta}(K)$, there exists $\overline{\gamma}\in \P^{top}_{\alpha,\beta}|K|^t_{reg}$ such that $\gamma = \overline{\gamma}\mu_{\ell}$. We obtain $\gamma\phi = \overline{\gamma}\mu_{\ell}\phi \mu^{-1}_{\ell}\mu_{\ell}$. Since $\mu_{\ell}\phi \mu^{-1}_{\ell}\in \mathcal{G}(1,1)$, we deduce that $\overline{\gamma}\mu_{\ell}\phi \mu^{-1}_{\ell}\in \P^{top}_{\alpha,\beta}|K|^t_{reg}$ and that $\gamma\phi\in \vec{R}^\ell_{\alpha,\beta}(K)$. 
\epf

\begin{nota}
	Let $\underline{x}=(x_1,\dots,x_n)$ and $\underline{x}'=(x'_1,\dots,x'_n)$ be two elements of $[0,1]^n$ with $n\geq 1$. Let
	\[
	d_1(\underline{x},\underline{x}') = \sum_{i=1}^{n} |x_i-x'_i|.
	\]
\end{nota}

An important feature shared by all $d$-paths (regular or not, tame or not) in the geometric realization of a precubical set is that they have a well-defined \textit{$L_1$-arc length} \cite[Section~2.2.1]{MR2521708} \cite[Section~2.2]{Raussen2012}. It is defined as follows. Consider a $d$-path $\gamma:[0,\ell]\to [0,1]^n$. Define the $L_1$-arc length between $\gamma(t)$ and $\gamma(t')$ by $d_1(\gamma(t),\gamma(t'))$. The $L_1$-arc length of a $d$-path of $[0,1]^n$ between $0_n$ an $1_n$ is therefore $n$. The $L_1$-arc length of a Moore composition of $d$-paths is defined by adding the $L_1$-arc length of each $d$-path. This definition makes sense since the coface operators preserve the metric $d_1$.

\begin{rem}
	The length of a nonconstant $d$-path $\gamma:[0,\ell]\to [0,1]^n$ between two vertices of $\{0,1\}^n$ is $\ell$, whereas its $L_1$-arc length is an integer belonging to $\{1,\dots,n\}$.
\end{rem} 

\bp \label{same-length}
Let $K$ be a precubical set. Let $(\alpha,\beta)\in K_0\p K_0$. Two execution paths of $\P^{top}_{\alpha,\beta}|K|^t_{reg}$ which are in the same path-connected component have the same $L_1$-arc length. 
\ep

\bpf
Let $\vec{P}_{\alpha,\beta}(K)_{co}$ be the set $\vec{P}_{\alpha,\beta}(K)$ equipped with the compact-open topology associated with the topology of $|K|_{geom}$. Write $(\P^{top}_{\alpha,\beta}|K|^t_{reg})_{co}$ for the underlying set of $\P^{top}_{\alpha,\beta}|K|^t_{reg}$ endowed with the compact-open topology. All elements of $(\P^{top}_{\alpha,\beta}|K|^t_{reg})_{co}$ are tame regular $d$-paths. We obtain a one-to-one continuous map \[(\P^{top}_{\alpha,\beta}|K|^t_{reg})_{co}\subset \vec{P}_{\alpha,\beta}(K)_{co}.\] By \cite[Proposition~2.2]{Raussen2012}, a composite map of the form \[[0,1]\longrightarrow (\P^{top}_{\alpha,\beta}|K|^t_{reg})_{co} \subset \vec{P}_{\alpha,\beta}(K)_{co} \longrightarrow \mathbb{R}\] is constant, where the right-hand map is the $L_1$-arc length. The $\Delta$-kelleyfication functor does not change the path-connected components. Hence the proof is complete. 
\epf

Due to the variety of the terminologies used in the mathematical literature, we recall the following definition for the convenience of the reader.

\bd
A \textit{pseudometric space} $(X,d)$ is a set $X$ equipped with a map $d:X\p X\to [0,+\infty]$ called a \textit{pseudometric} such that:
\begin{itemize}
	\item $\forall x\in X,d(x,x)=0$
	\item $\forall (x,y)\in X\p X, d(x,y)=d(y,x)$
	\item $\forall (x,y,z)\in X\p X\p X, d(x,y)\leq d(x,z)+d(z,y)$.
\end{itemize}
A map $f:(X,d)\to (Y,d)$ of pseudometric spaces is a set map $f:X\to Y$ which does not increase distance, i.e. $\forall (x,y)\in X\p X, d(f(x),f(y))\leq d(x,y)$.
\ed

Every metric space is a pseudometric space. The interest of this category for Proposition~\ref{cont} lies in the following two facts. At first, it is bicomplete, being a coreflective subcategory of the bicomplete category of Lawvere metric spaces which are pseudometric spaces with the symmetry condition dropped \cite{LawvereMetric}. Secondly, the family of balls $B(x,\epsilon)=\{y\in X\mid d(x,y)<\epsilon\})$ with $x\in X$ and $\epsilon>0$ generates a topology called the underlying topology of $(X,d)$. This gives rise to a functor from pseudometric spaces to general topological spaces.

By equipping each topological $n$-cube $[0,1]^n$ for $n\geq 0$ with the metric $d_1$, we obtain a cocubical pseudometric space, and by Proposition~\ref{eq-cat}, a colimit-preserving functor from precubical sets to pseudometric spaces. By composing with the underlying topological space functor, we obtain a functor from precubical sets to general topological spaces denoted by $K\mapsto |K|_{d_1}$. Since there is a natural homeomorphism $|\square[n]|_{geom} \iso |\square[n]|_{d_1}$ for all $n\geq 0$, the universal property of the colimit proves that the identity induces a continuous map \[|K|_{geom} \longrightarrow |K|_{d_1}.\] It is worth noting that the latter map is a homeomorphism if and only if the CW-complex $|K|_{geom}$ is locally finite. Indeed, it is a homeomorphism if and only if the final topology is pseudometrizable, and therefore, if and only if each path-connected component is metrizable. We conclude using \cite[Proposition~1.5.17]{MR1074175}. Note that this argument also implicitly proves that the topology of $|K|_{geom}$ is pseudometrizable if and only if it is metrizable (in this case, there exists a real number $B>0$ independent from the path-connected components such that each path-connected component of $|K|_{geom}$ is included in a ball of radius $B$, being embeddable in the Hilbert cube: see \cite[Proposition~1.5.12, Proposition~1.5.13 and Theorem~1.5.16]{MR1074175}).

\bp \label{cont}
Let $K$ be a precubical set. Let $\ell >0$. The function \[L:\vec{R}^\ell(K) \p [0,\ell] \longrightarrow [0,+\infty[\] which takes $(\gamma,t)$ to the $L_1$-arc length between $\gamma(0)$ and $\gamma(t)$ is continuous.
\ep

\bpf
Let $\vec{P}_{\alpha,\beta}(K)_{d_1}$ be the set $\vec{P}_{\alpha,\beta}(K)$ equipped with the compact-open topology associated with the topological space $|K|_{d_1}$. Using \cite[Lemma~2.13]{MR2521708}, the set map $L:\vec{P}(K)_{d_1} \p [0,+\infty[ \longrightarrow [0,+\infty[$ taking a pair $(\gamma,t)$ to the $L_1$-arc length between $\gamma(0)$ and $\gamma(t)$ is continuous. Let $(\vec{R}^\ell(K))_{d_1}$ ($(\vec{R}^\ell(K))_{co}$ resp.) be the underlying set of the space $\vec{R}^\ell(K)$ equipped with the compact-open topology associated with the topological space $|K|_{d_1}$ (associated with the topological space $|K|_{geom}$ resp.). Since the identity induces a continuous map $|K|_{geom} \to |K|_{d_1}$, we obtain a continuous map 
\[L:(\vec{R}^\ell(K))_{co} \p [0,\ell] \to (\vec{R}^\ell(K))_{d_1} \p [0,\ell] \subset \vec{P}(K)_{d_1} \p [0,+\infty[\to [0,+\infty[.\]
Finally, take the image by the $\Delta$-kelleyfication functor. The latter is a right adjoint therefore it preserves binary products. Besides, $[0,\ell]$ and $[0,+\infty[$ are already $\Delta$-generated. Hence the proof is complete.
\epf

\begin{nota}
	By adjunction, we obtain a continuous map \[\widehat{L}:\vec{R}^\ell(K) \longrightarrow \ttop([0,\ell],[0,+\infty[).\] 
\end{nota}

\begin{lem} \label{minicalcul}
	For all $r\in \vec{R}^\ell(K)$, for all $\phi\in \mathcal{G}(\ell',\ell)$, and for all $t\in [0,\ell']$, one has \[\widehat{L}(r\phi)(t) = \widehat{L}(r)(\phi(t)).\]
\end{lem}

\bpf
The $L_1$-arc length between $r(0)=r(\phi(0))$ and $r(\phi(t))$ for the $d$-path $r$ is equal to the $L_1$-arc length between $(r\phi)(0)$ and $(r\phi)(t)$ for the $d$-path $r\phi$.
\epf

Intuitively, the natural $d$-paths are the $d$-paths whose speed corresponds to the $L_1$-arc length.

\bd \cite[Definition~2.14]{MR2521708}
Let $K$ be a precubical set. Let $(\alpha,\beta)\in K_0\p K_0$. A $d$-path $\gamma$ of $\vec{R}^\ell_{\alpha,\beta}(K)$ is \textit{natural} if $\widehat{L}(\gamma)(t)=t$ for all $t\in [0,\ell]$. This implies that $\ell$ is an integer (greater than or equal to $1$). The subset of natural $d$-paths of length $n\geq 1$ from $\alpha$ to $\beta$ equipped with the $\Delta$-kelleyfication of the relative topology is denoted by $\vec{N}^n_{\alpha,\beta}(K)$. Let 
\[
\vec{N}^\ell(K) = \coprod_{(\alpha,\beta)\in K_0\p K_0} \vec{N}^\ell_{\alpha,\beta}(K).
\]
\ed

The following theorem is an improvement and an adaptation in our topological setting of \cite[Proposition~2.16]{MR2521708}: the homotopy equivalence is replaced by a homeomorphism thanks to the spaces $\mathcal{G}(\ell,n)$. 

\bth \label{Psi}
Let $K$ be a precubical set. Let $(\alpha,\beta)\in K_0\p K_0$. There is a homeomorphism 
\[
\Psi^{\ell}:\vec{R}_{\alpha,\beta}^{\ell}(K) \stackrel{\iso}\longrightarrow \displaystyle\coprod_{n\geq 1} \mathcal{G}(\ell,n) \p \vec{N}_{\alpha,\beta}^n(K).
\]
\eth

\bpf
By Proposition~\ref{same-length}, the space $\vec{R}_{\alpha,\beta}^{\ell}(K)$ is the direct sum of the subspaces $\vec{R}_{\alpha,\beta}^{\ell,n}(K)$ of tame regular $d$-paths of $L_1$-arc length $n$ for $n\geq 1$. It then suffices to prove the homeomorphism $\vec{R}_{\alpha,\beta}^{\ell,n}(K) \iso \mathcal{G}(\ell,n) \p \vec{N}_{\alpha,\beta}^n(K)$ for all $n\geq 1$. Consider the continuous map $\Phi^{\ell}:\mathcal{G}(\ell,n)\p \vec{N}_{\alpha,\beta}^n(K) \to \vec{R}_{\alpha,\beta}^{\ell,n}(K)$ defined by \[\Phi^{\ell}(\phi,\gamma)=\gamma\phi.\] Since $\widehat{L}(\gamma\phi)=\widehat{L}(\gamma)\phi=\phi$, the first equality by Lemma~\ref{minicalcul} and the second equality since $\gamma$ is natural, the $d$-path $\gamma\phi$ has no stop-interval, $\phi$ being bijective, and it is therefore regular. We want to define a continuous map $\Psi^{\ell}:\vec{R}_{\alpha,\beta}^{\ell,n}(K)\to \mathcal{G}(\ell,n)\p \vec{N}_{\alpha,\beta}^n(K)$. Let $r\in \vec{R}_{\alpha,\beta}^{\ell,n}(K)$. Then $\widehat{L}(r):[0,\ell] \to [0,n]$ is surjective and nondecreasing. Let $t,t'\in [0,1]$ such that $\widehat{L}(r)(t)=\widehat{L}(r)(t')$. Since $r$ is regular, it has no stop interval. It implies that $t=t'$. Thus, $\widehat{L}(r)\in \mathcal{G}(\ell,n)$. Let \[\Psi^{\ell}(r)=(\widehat{L}(r),r\widehat{L}(r)^{-1}).\] The map $\Psi^{\ell}$ is continuous by Proposition~\ref{cont} and by \cite[Lemma~6.2]{Moore2}. Then one has \[\widehat{L}(r\widehat{L}(r)^{-1})(t)=\widehat{L}(r)(\widehat{L}(r)^{-1}(t)) = t\] for all $0\leq t \leq n$, the first equality by Lemma~\ref{minicalcul} and the second equality by algebraic simplification. This means that $r\widehat{L}(r)^{-1}\in \vec{N}_{\alpha,\beta}^n(K)$: the path $r\widehat{L}(r)^{-1}$ is called the naturalization of $r$ in \cite[Definition~2.14]{MR2521708}. One has \[\Phi^{\ell}\Psi^{\ell}(r)=\Phi^{\ell}(\widehat{L}(r),r\widehat{L}(r)^{-1})= r\widehat{L}(r)^{-1}\widehat{L}(r) = r,\] the first equality by definition of $\Psi^{\ell}$, the second equality by definition of $\Phi^{\ell}$ and the last equality by algebraic simplification. Let $(\phi,\gamma)\in \mathcal{G}(\ell,n) \p \vec{N}_{\alpha,\beta}^n(K)$. Since $\widehat{L}(\gamma\phi) = \phi$, we obtain \[\Psi^{\ell}\Phi^{\ell}(\phi,\gamma)=\Psi^{\ell}(\gamma\phi) = (\widehat{L}(\gamma\phi),\gamma\phi\widehat{L}(\gamma\phi)^{-1}) = (\phi,\gamma),\]
the first equality by definition of $\Phi^{\ell}$, the second equality by definition of $\Psi^{\ell}$, and the last equality by algebraic simplification.
\epf

Theorem~\ref{Psi} can be applied as follows.

\bd \label{rea_conc} \cite[Definition~7.1]{NaturalRealization}
Let $K$ be a precubical set. The \textit{tame concrete realization} of $K$ is the flow $|K|_{tc}$ defined as follows:
\begin{align*}
	& |K|_{tc}^0=K_0 \\
	& \forall (\alpha,\beta)\in K_0\p K_0, \P_{\alpha,\beta} |K|_{tc}=\coprod_{n\geq 1}\vec{N}_{\alpha,\beta}^n(K)
\end{align*}
The composition of execution paths is induced by the Moore composition.
\ed

As noticed after \cite[Theorem~7.7]{NaturalRealization}, the tame concrete realization functor is not colimit-preserving in general. It is therefore not a realization functor for flows in the sense of \cite[Definition~3.6]{NaturalRealization}. This means that it is difficult to calculate for a general precubical set. However, \cite[Theorem~7.7]{NaturalRealization} proves that it coincides with a colimit-realization functor when the precubical set is spatial in the sense of Definition~\ref{carac_spatial}. Since most of the concrete examples coming from concurrency theory are spatial precubical sets (see the comment after Definition~\ref{carac_spatial}), the tame concrete realization of a precubical set is easily calculable for this kind of examples.

\bp \label{same}
For every precubical set $K$, there is the natural isomorphism of flows \[\Cat(|K|^t_{reg})\iso |K|_{tc}.\]
\ep

\bpf 
By Theorem~\ref{Psi}, there is the homeomorphism
\[
\P^{top}_{\alpha,\beta}|K|^t_{reg} = \vec{R}_{\alpha,\beta}^{1}(K) \stackrel{\iso}\longrightarrow \displaystyle\coprod_{n\geq 1} \mathcal{G}(1,n) \p \vec{N}_{\alpha,\beta}^n(K)
\]
for all $(\alpha,\beta)\in K_0\p K_0$. The coequalizer above identifies $(\psi,\gamma)\in \mathcal{G}(1,n) \p \vec{N}_{\alpha,\beta}^n(K)$ with $(\psi\phi,\gamma)\in \mathcal{G}(1,n) \p \vec{N}_{\alpha,\beta}^n(K)$ for all $\phi\in \mathcal{G}(1,1)$. The proof is complete by Definition~\ref{rea_conc}.
\epf

\section{Cube chains and tame regular realization}
\label{cubechain}

\bd
The cocubical Moore flow $\moore^{\mathcal{G}}|\square[*]|^t_{reg}$ gives rise by Proposition~\ref{eq-cat} to a colimit-preserving functor $[-]_{reg} : \square^{op}\set\to \dtopG$ defined by 
\[
[K]_{reg} = \liminj_{\square[n]\rightarrow K} \moore^{\mathcal{G}}|\square[n]|^t_{reg}
\]
\ed

The following notations coincide with \cite[Definition~4.8 and Definition~4.12]{NaturalRealization} by \cite[Proposition~4.10 and Proposition~4.14]{NaturalRealization}.

\begin{nota}
	Let $N_n = \vec{N}_{0_n,1_n}^n(\square[n])$ and $\de N_n = \vec{N}_{0_n,1_n}^n(\de\square[n])$ for $n\geq 0$. These spaces are first countable, $\Delta$-generated and $\Delta$-Hausdorff.
\end{nota}

The continuous map \[[0,1]^{m_1}\sqcup [0,1]^{m_2} \longrightarrow [0,1]^{m_1+m_2}\] defined by taking the tuple $(t_1,\dots,t_{m_1})$ to the tuple $(t_1,\dots,t_{m_1},0_{m_2})$ and the tuple $(t'_1,\dots,t'_{m_2})$ to the tuple $(1_{m_1},t'_1,\dots,t'_{m_2})$ induces a continuous map \[N_{m_1}\p N_{m_2}\longrightarrow N_{m_1+m_2}\] by using the fact that the Moore composition of two natural $d$-paths is still a natural $d$-path.

Cube chains are a powerful notion introduced by Ziemia{\'{n}}ski. We use the presentation given in \cite[Section~7]{MR4070250}. 

\begin{nota}
	Let $\seq(n)$ be the set of sequences of positive integers \[\underline{n}=(n_1,\dots,n_p)\] with $n_1+\dots + n_p=n$. Let $\underline{n}=(n_1,\dots,n_p) \in \seq(n)$. Then $|\underline{n}|=n$ is the \textit{length} of $\underline{n}$ and $\ell(\underline{n})=p$ is the \textit{number of elements} of $\underline{n}$. 
\end{nota}

Let $K$ be a precubical set and $A={a_1<\dots<a_k} \subset \{1,\dots,n\}$ and $\epsilon\in \{0,1\}$. The \textit{iterated face map} is defined by \[\de^\epsilon_A=\de^\epsilon_{a_1} \de^\epsilon_{a_2} \dots \de^\epsilon_{a_k}.\]

\bd 
Let $\underline{n}=(n_1,\dots,n_p)\in \seq(n)$. The \textit{$\underline{n}$-cube} is the precubical set 
\[
\square[\underline{n}] = \square[n_1] * \dots * \square[n_p]
\]
where the notation $*$ means that the final state $1_{n_i}$ of the precubical set $\square[n_i]$ is identified with the initial state $0_{n_{i+1}}$ of the precubical set $\square[n_{i+1}]$ for $1\leq i\leq p-1$.
\ed 

Let $K$ be a precubical set. Let $\alpha,\beta\in K_0$. Let $n\geq 1$. The small category $\Ch_{\alpha,\beta}(K,n)$ is defined as follows. The objects are the maps of precubical sets $\square[\underline{n}] \to K$ with $|\underline{n}|=n$ where the initial state of $\square[n_1]$ is mapped to $\alpha$ and the final state of $\square[n_p]$ is mapped to $\beta$. Let $A\sqcup B=\{1,\dots,m_1+m_2\}$ be a partition with the cardinal of $A$ equal to $m_1>0$ and the cardinal of $B$ equal to $m_2>0$. Let \[\phi_{A,B}:\square[m_1]*\square[m_2] \longrightarrow \square[m_1+m_2]\] be the unique map of precubical sets such that 
\begin{align*}
	&\phi_{A,B}(\id_{[m_1]}) = \de^0_B(\id_{[m_1+m_2]}),\\
	&\phi_{A,B}(\id_{[m_2]}) = \de^1_A(\id_{[m_1+m_2]}).
\end{align*}
For $i\in \{1,\dots,\ell(\underline{n})\}$ and a partition $A\sqcup B=\{1,\dots,n_i\}$, let \[\delta_{i,A,B}=\id_{\square[n_1]}*\dots*\id_{\square[n_{i-1}]}*\phi_{A,B}*\id_{\square[n_{i+1}]}*\dots*\id_{\square[n_{\ell(\underline{n})}]}.\] The morphisms are the commutative diagrams
\[
\xymatrix@C=3em@R=3em
{
	\square[\underline{n}_a] \fd{}\fr{a} & K \ar@{=}[d] \\
	\square[\underline{n}_b]  \fr{b}  & K
}
\] 
where the left vertical map is a composite of maps of precubical sets of the form $\delta_{i,A,B}$. 

We recall the definition of the functor $K\mapsto ||K||$ from precubical sets to flows introduced in \cite[Section~6]{NaturalRealization}. The set of states of $||K||$ is $K_0$. Consider the small diagram of spaces \[\D_{\alpha,\beta}(K,n):\Ch_{\alpha,\beta}(K,n)\longrightarrow \top\] defined by on objects by \[\D_{\alpha,\beta}(K,n)(\square[\underline{n}]\to K) = N_{n_1} \p \dots \p N_{n_p}\] with $\underline{n} = (n_1,\dots,n_p)$ and $\sum_i n_i=n$ and on morphisms by using the maps \[\P|\phi_{A,B}|_{nat}:\P(\square[m_1]*\square[m_2]) \to \P(\square[m_1+m_2])\] which induce maps $N_{m_1}\p N_{m_2}\to N_{m_1+m_2}$ given by the Moore composition of tame natural $d$-paths. The space of execution spaces $\P_{\alpha,\beta} ||K||$ is defined as follows:
\[
\P_{\alpha,\beta} ||K|| = \displaystyle\coprod_{n\geq 1} \liminj \D_{\alpha,\beta}(K,n).
\]
The concatenation of tuples induces functors \[\D_{\alpha,\beta}(K,m_1)\p \D_{\beta,\gamma}(K,m_2) \to \D_{\alpha,\gamma}(K,m_1+m_2),\]
and, using \cite[Proposition~A.4]{leftproperflow}, continuous maps 
\[\liminj \D_{\alpha,\beta}(K,m_1)\p \liminj \D_{\beta,\gamma}(K,m_2) \to \liminj\D_{\alpha,\gamma}(K,m_1+m_2)\] for all $m_1,m_2\geq 1$. We obtain an associative composition map \[\P_{\alpha,\beta} ||K|| \p \P_{\beta,\gamma} ||K|| \to \P_{\alpha,\gamma} ||K||\] for all $(\alpha,\beta,\gamma)\in K_0\p K_0\p K_0$.

We want to define a Moore flow $||K||^{\mathcal{G}}$ by mimicking the above construction of $||K||$ by using the following two rules: 1) any occurrence of the topological space $N_k$ is replaced by the $\mathcal{G}$-space $\mathbb{F}^{\mathcal{G}^{op}}_k N_k$ for all integers $k\geq 1$; 2) any product of spaces of the form $N_{n_1}\p\dots \p N_{n_p}$ is replaced by the tensor products of $\mathcal{G}$-spaces  $\mathbb{F}^{\mathcal{G}^{op}}_{n_1}N_{n_1}\ot\dots \ot \mathbb{F}^{\mathcal{G}^{op}}_{n_p}N_{n_p}$. The idea is to work with the biclosed semimonoidal category $(\topdgr_0,\ot)$ instead of with the biclosed (semi)monoidal category $(\top,\p)$.

Let $(\alpha,\beta)\in K_0\p K_0$. Let \[\D_{\alpha,\beta}^{\mathcal{G}}(K,n):\Ch_{\alpha,\beta}(K,n)\longrightarrow \topdgr_0\] be the functor defined on objects by \[\D_{\alpha,\beta}^{\mathcal{G}}(K,n)(\square[\underline{n}]\to K) = \mathbb{F}^{\mathcal{G}^{op}}_{n_1} N_{n_1} \ot \dots \ot \mathbb{F}^{\mathcal{G}^{op}}_{n_p} N_{n_p}\] with $\underline{n} = (n_1,\dots,n_p)$ and $\sum_i n_i=n$ and on morphisms by taking the map \[\phi_{A,B}:\square[m_1]*\square[m_2] \to \square[m_1+m_2]\] to the composite map of $\mathcal{G}$-spaces \[\mathbb{F}^{\mathcal{G}^{op}}_{m_1} N_{m_1}\ot \mathbb{F}^{\mathcal{G}^{op}}_{m_1} N_{m_2} \iso \mathbb{F}^{\mathcal{G}^{op}}_{m_1+m_2} (N_{m_1}\p N_{m_2}) \to \mathbb{F}^{\mathcal{G}^{op}}_{m_1+m_2} N_{m_1+m_2}\] where the isomorphism is given by Proposition~\ref{Ftenseur} and where the map $N_{m_1}\p N_{m_2} \to N_{m_1+m_2}$ is given by the Moore composition of tame natural $d$-paths. The $\mathcal{G}$-space of execution spaces $\P_{\alpha,\beta} ||K||^{\mathcal{G}}$ is defined as follows:
\[
\P_{\alpha,\beta} ||K||^{\mathcal{G}} = \displaystyle\coprod_{n\geq 1} \liminj \D_{\alpha,\beta}^{\mathcal{G}}(K,n).
\]
There is the obvious proposition: 
\bp
Let $D_i:I_i\to \topdgr_0$ be two small diagrams of $\mathcal{G}$-spaces with $i=1,2$. Then the mappings ($\Mor(I)$ means the set of morphisms of a small category $I$)
\begin{align*}
	& (i_1,i_2)\in I_1\p I_2 \mapsto D_1(i_1)\ot D_2(i_2)\\
	& (f,g)\in \Mor(I_1)\p \Mor(I_2)\mapsto f\ot g
\end{align*} yield a well defined small diagram of $\mathcal{G}$-spaces denoted by \[D_1\ot D_2:I_1\p I_2 \longrightarrow \topdgr_0.\] 
\ep

Proposition~\ref{cassepied} is the main fact which enables us to define the composition law of the Moore flow $||K||^{\mathcal{G}}$.

\bp \label{cassepied}
Let $K$ be a precubical set. Let $\alpha,\beta,\gamma\in K_0$. Let $m_1,m_2\geq 1$. There is the isomorphism of $\mathcal{G}$-spaces 
\[
\liminj \big(\D_{\alpha,\beta}^{\mathcal{G}}(K,m_1)\ot \D_{\beta,\gamma}^{\mathcal{G}}(K,m_2)\big) \iso \big(\liminj\D_{\alpha,\beta}^{\mathcal{G}}(K,m_1)\big)\ot \big(\liminj\D_{\beta,\gamma}^{\mathcal{G}}(K,m_2)\big).
\]
\ep

\bpf
The semimonoidal category $(\topdgr_0,\ot)$ being biclosed by Theorem~\ref{closedsemimonoidal}, write
\begin{align}
	&\topdgr_0(D,\{E,F\}_L) \iso \topdgr_0(D\ot E,F), \label{left} \tag{L}\\
	&\topdgr_0(E,\{D,F\}_R) \iso \topdgr_0(D\ot E,F), \label{right} \tag{R}
\end{align}
$D,E,F$ being three $\mathcal{G}$-spaces. Then we obtain the sequence of bijections
\begin{align*}
	\topdgr_0&\big(\big(\liminj\D_{\alpha,\beta}^{\mathcal{G}}(K,m_1)\big)\ot \big(\liminj\D_{\beta,\gamma}^{\mathcal{G}}(K,m_2)\big),F\big) \\
	&\iso \topdgr_0\big(\liminj\D_{\alpha,\beta}^{\mathcal{G}}(K,m_1),\big\{\liminj \D_{\beta,\gamma}^{\mathcal{G}}(K,m_2),F\big\}_L\big)\\
	&\iso \limproj_{a_1} \topdgr_0\big(\D_{\alpha,\beta}^{\mathcal{G}}(K,m_1)(a_1),\big\{\liminj \D_{\beta,\gamma}^{\mathcal{G}}(K,m_2),F\big\}_L\big)\\
	&\iso \limproj_{a_1} \topdgr_0\big(\D_{\alpha,\beta}^{\mathcal{G}}(K,m_1)(a_1)\ot \big(\liminj\D_{\beta,\gamma}^{\mathcal{G}}(K,m_2)\big),F\big)\\
	&\iso \limproj_{a_1} \topdgr_0\big(\liminj \D_{\beta,\gamma}^{\mathcal{G}}(K,m_2),\big\{\D_{\alpha,\beta}^{\mathcal{G}}(K,m_1)(a_1),F\big\}_R\big) \\
	&\iso \limproj_{a_1} \limproj_{a_2} \topdgr_0\big( \D_{\beta,\gamma}^{\mathcal{G}}(K,m_2)(a_2),\big\{\D_{\alpha,\beta}^{\mathcal{G}}(K,m_1)(a_1),F\big\}_R\big) \\
	&\iso \limproj_{a_1} \limproj_{a_2} \topdgr_0\big(\D_{\alpha,\beta}^{\mathcal{G}}(K,m_1)(a_1)\ot \D_{\beta,\gamma}^{\mathcal{G}}(K,m_2)(a_2),F\big)\\
	&\iso \limproj_{(a_1,a_2)} \topdgr_0\big(\D_{\alpha,\beta}^{\mathcal{G}}(K,m_1)(a_1)\ot \D_{\beta,\gamma}^{\mathcal{G}}(K,m_2)(a_2),F\big)\\
	&\iso \topdgr_0\big(\liminj \big(\D_{\alpha,\beta}^{\mathcal{G}}(K,m_1)\ot \D_{\beta,\gamma}^{\mathcal{G}}(K,m_2)\big),F\big)
\end{align*}
for all $\mathcal{G}$-spaces $F$, the first and third bijections by \eqref{left}, the second and fifth and eighth bijections by the universal property of the colimit, the fourth and sixth bijections by \eqref{right} and finally the seventh bijection because limits commute with each other. The proof is complete thanks to the Yoneda lemma. 
\epf

We then consider the category of \textit{all} small diagrams of $\mathcal{G}$-spaces over \textit{all} small categories, denoted by $\mathrm{Diag}(\topdgr_0)$, defined as follows. An object is a functor $F:{I}\to \topdgr_0$. A morphism from $F:{I}_1\to \topdgr_0$ to $G:{I}_2\to \topdgr_0$ is a pair $(f:{I}_1\to {I}_2,\mu:F \Rightarrow G.f)$ where $f$ is a functor and $\mu$ is a natural transformation. If $(g,\nu)$ is a map from $G:{I}_2\to \topdgr_0$ to $H:\underline{K}\to \topdgr_0$, then the composite $(g,\nu).(f,\mu)$ is defined by $(g.f,(\nu.f)\odot\mu)$, $\odot$ meaning the composition of natural transformations. The identity of $F:{I}_1\to \topdgr_0$ is the pair $(\id_{{I}_1},\id_F)$. It is well-known that this defines an associative composition law (e.g. see \cite[Appendix~A]{leftproperflow}). The colimit construction induces a functor \[\liminj :\mathrm{Diag}(\topdgr_0) \longrightarrow \topdgr_0\] by \cite[Proposition~A.2]{leftproperflow}. We define a map of $\mathrm{Diag}(\topdgr_0)$
\[(f^{m_1,m_2}_{\alpha,\beta,\gamma},\mu^{m_1,m_2}_{\alpha,\beta,\gamma}):\D_{\alpha,\beta}^{\mathcal{G}}(K,m_1)\ot \D_{\beta,\gamma}^{\mathcal{G}}(K,m_2) \longrightarrow \D_{\alpha,\gamma}^{\mathcal{G}}(K,m_1+m_2)\]
for all $m_1,m_2\geq 1$ and all $(\alpha,\beta,\gamma)\in K_0\p K_0\p K_0$ as follows. The functor \[f^{m_1,m_2}_{\alpha,\beta,\gamma}:\Ch_{\alpha,\beta}(K,m_1)\p \Ch_{\beta,\gamma}(K,m_2) \to \Ch_{\alpha,\gamma}(K,m_1+m_2)\] takes a pair $(\square[\underline{m_1}]\to K,\square[\underline{m_2}] \to K)$ to the map of precubical sets $\square[\underline{m_1}]*\square[\underline{m_2}]\to K$ and the natural transformation
\[
\mu^{m_1,m_2}_{\alpha,\beta,\gamma}:\D_{\alpha,\beta}^{\mathcal{G}}(K,m_1)\ot \D_{\beta,\gamma}^{\mathcal{G}}(K,m_2)\Longrightarrow \D_{\alpha,\gamma}^{\mathcal{G}}(K,m_1+m_2) f^{m_1,m_2}_{\alpha,\beta,\gamma}
\] 
is the identity. Using Proposition~\ref{cassepied}, we obtain a map of $\mathcal{G}$-spaces \[\P_{\alpha,\beta} ||K||^{\mathcal{G}} \ot \P_{\beta,\gamma} ||K||^{\mathcal{G}} \longrightarrow \P_{\alpha,\gamma} ||K||^{\mathcal{G}}\] for all $(\alpha,\beta,\gamma)\in K_0\p K_0\p K_0$. It is strictly associative since $(\topdgr_0,\ot)$ is semimonoidal \cite[Proposition~5.11]{Moore1}. We have obtained a well-defined Moore flow $||K||^{\mathcal{G}}$.

\bp \label{cp}
The mapping $(\phi\gamma)\mapsto \gamma\phi$ yields the homeomorphisms
\begin{align*}
	&\mathcal{G}(\ell,n)\p N_n \iso \P_{0_{n},1_{n}}^\ell\moore^{\mathcal{G}}|\square[n]|^t_{reg} \\
	&\mathcal{G}(\ell,n)\p \de N_n \iso \P_{0_{n},1_{n}}^\ell\moore^{\mathcal{G}}|\de\square[n]|^t_{reg}
\end{align*}
for all $n\geq 1$ and all $\ell>0$. In particular, for $\ell=1$, we obtain 
the homeomorphisms
\begin{align*}
	&\mathcal{G}(1,n)\p N_n \iso \P^{top}_{0_{n},1_{n}}|\square[n]|^t_{reg} \\
	&\mathcal{G}(1,n)\p \de N_n \iso \P^{top}_{0_{n},1_{n}}|\de\square[n]|^t_{reg}
\end{align*}
for all $n\geq 1$. Note that for $n=1$, $\de N_1=\varnothing= \P^{top}_{0_{n},1_{n}}|\de\square[n]|^t_{reg}$.
\ep

\bpf
The first part is a consequence of Theorem~\ref{Psi} applied to the regular $d$-paths of length $\ell$ and of $L_1$-arc length $n$ of the precubical sets $\square[n]$ and $\de\square[n]$. The particular case $\ell=1$ is a consequence of the definition of the functor $\moore^{\mathcal{G}}:\ptop{\mathcal{G}}\to\dtopG$
\epf

\begin{cor} \label{cp2}
	For all $n\geq 1$, there are the isomorphisms of $\mathcal{G}$-spaces 
	\begin{align*}
		&\mathbb{F}^{\mathcal{G}^{op}}_nN_n\iso \P_{0_{n},1_{n}}\moore^{\mathcal{G}}|\square[n]|^t_{reg},\\
		&\mathbb{F}^{\mathcal{G}^{op}}_n \de N_n\iso \P_{0_{n},1_{n}}\moore^{\mathcal{G}}|\de\square[n]|^t_{reg}.
	\end{align*}
\end{cor}

\bpf
It is a consequence of Proposition~\ref{cp}, of the fact that all regular $d$-paths from $0_{n}$ to $1_{n}$ of $|\square[n]|$ are of $L_1$-arc length $n$, and of the definition of the functor $\moore^{\mathcal{G}}:\ptop{\mathcal{G}}\to\dtopG$. 
\epf

\bth \label{doublebarG}
There is a natural isomorphism of Moore flows \[[K]_{reg} \iso ||K||^{\mathcal{G}}\] for all precubical sets $K$.
\eth

The proof of Theorem~\ref{doublebarG} is roughly speaking the proof of \cite[Proposition~6.2]{NaturalRealization} and \cite[Theorem~6.3]{NaturalRealization} by working with the biclosed semimonoidal category $(\topdgr_0,\ot)$ instead of with the biclosed (semi)monoidal category $(\top,\p)$.  Some details slightly change due to the fact that we work with Moore composition of paths which have a length.

\bpf
First of all, we prove that there is an isomorphism of cocubical Moore flows $[\square[*]]_{reg} = \moore^{\mathcal{G}}|\square[*]|^t_{reg} \iso ||\square[*]||^{\mathcal{G}}$. Let $d_1(\alpha,\beta)=m$. For $\alpha<\beta\in \{0,1\}^n$ for the product order, the small category $\Ch_{\alpha,\beta}(\square[n],p)$ is empty if $p\neq m$ and it has a terminal object $\square[m]\to \square[n]$ corresponding to the subcube from $\alpha$ to $\beta$ of $\square[n]$ if $p= m$. We deduce the isomorphisms of $\mathcal{G}$-spaces 
\begin{align*}
	\P_{\alpha,\beta} ||\square[n]||^{\mathcal{G}} &= \liminj_{\substack{\underline{n}=(n_1,\dots,n_p),\ell(\underline{n})=m\\\square[\underline{n}]\to \square[n] \in \Ch_{\alpha,\beta}(\square[n],m)}} \mathbb{F}^{\mathcal{G}^{op}}_{n_1}N_{n_1} \ot \dots \ot \mathbb{F}^{\mathcal{G}^{op}}_{n_p}N_{n_p} \\
	&\iso \liminj_{\substack{\underline{n}=(n_1,\dots,n_p),\ell(\underline{n})=m\\\square[\underline{n}]\to \square[n] \in \Ch_{\alpha,\beta}(\square[n],m)}} \mathbb{F}^{\mathcal{G}^{op}}_{m}(N_{n_1} \p \dots \p N_{n_p})\\
	&\iso \mathbb{F}^{\mathcal{G}^{op}}_{m}\bigg(\liminj_{\substack{\underline{n}=(n_1,\dots,n_p),\ell(\underline{n})=m\\\square[\underline{n}]\to \square[n] \in \Ch_{\alpha,\beta}(\square[n],m)}} (N_{n_1} \p \dots \p N_{n_p})\bigg)\\
	&\iso \mathbb{F}^{\mathcal{G}^{op}}_{m} N_m \\&\iso \P_{\alpha,\beta} [\square[n]]_{reg},
\end{align*}
the first equality by definition of $||\square[n]||^{\mathcal{G}}$, the first isomorphism by Proposition~\ref{Ftenseur}, the second isomorphism by Proposition~\ref{ev-adj}, the third isomorphism because of the unique map $\underline{c}:\square[m]\to \square[n]$ which is the terminal object of $\Ch_{\alpha,\beta}(\square[n],m)$, and the last isomorphism by Corollary~\ref{cp2} and by definition of the Moore flow $[\square[n]]_{reg}$. The universal property of the colimit yields a natural map of Moore flows $[K]_{reg} \to ||K||^{\mathcal{G}}$, the functor $[-]_{reg}$ being colimit-preserving. Let $\underline{n}=(n_1,\dots,n_p) \in \seq(n)$. Every map of precubical sets $\square[\underline{n}]\to K$ gives rise to a map of Moore flows $[\square[\underline{n}]]_{reg}\to [K]_{reg}$. Let $\alpha_0,\alpha_1,\dots,\alpha_p\in K_0$ the images by $0_{n_1}, 1_{n_1}=0_{n_2},\dots,1_{n_p}$ by this map of precubical sets. The are maps of Moore flows $[\square[n_i]]_{reg} \to [K]_{reg}$ inducing maps of $\mathcal{G}$-spaces \[\mathbb{F}^{\mathcal{G}^{op}}_{n_i}N_{n_i}\iso \P_{0_{n_i},1_{n_i}}[\square[n_i]]_{reg} \longrightarrow \P_{\alpha_{i-1},\alpha_i}[K]_{reg}\] for $i\in \{1,\dots,p\}$. Using the composition law of the Moore flow $[K]_{reg}$, one obtains a map of $\mathcal{G}$-spaces  \[\mathbb{F}^{\mathcal{G}^{op}}_{n_1}N_{n_1}\ot \dots \ot  \mathbb{F}^{\mathcal{G}^{op}}_{n_p}N_{n_p} \longrightarrow \P  [K]_{reg}.\]
Consequently, we obtain a cocone \[(\mathbb{F}^{\mathcal{G}^{op}}_{n_1}N_{n_1}\ot \dots \ot  \mathbb{F}^{\mathcal{G}^{op}}_{n_p}N_{n_p})_{\substack{\square[\underline{n}]\to K \\\in \Ch_{\alpha,\beta}(K,n)}} \stackrel{\bullet}\longrightarrow \P [K]_{reg}\] and then a map of Moore flows $||K||^{\mathcal{G}} \to [K]_{reg}$ which is bijective on states. The composite map of Moore flows $[K]_{reg}\to ||K||^{\mathcal{G}} \to [K]_{reg}$ is the identity of $[K]_{reg}$ because it is the identity for $K=\square[n]$ for all $n\geq 0$. Consequently, for all $(\alpha,\beta)\in K_0\p K_0$, the composite continuous map $\P_{\alpha,\beta}^1[K]_{reg}\to \P_{\alpha,\beta}^1||K||^{\mathcal{G}} \to \P_{\alpha,\beta}^1[K]_{reg}$ is the identity of $\P_{\alpha,\beta}^1[K]_{reg}$: this means that the left-hand map $\P_{\alpha,\beta}^1[K]_{reg}\to \P_{\alpha,\beta}^1||K||^{\mathcal{G}}$ is one-to-one and that the right-hand map $\P_{\alpha,\beta}^1||K||^{\mathcal{G}} \to \P_{\alpha,\beta}^1[K]_{reg}$ is onto. Consider an element $\gamma\in \P_{\alpha,\beta}^1||K||^{\mathcal{G}}$. It has a representative of the form $\overline{\gamma}\in \mathbb{F}^{\mathcal{G}^{op}}_{n_1}N_{n_1}\ot \dots \ot  \mathbb{F}^{\mathcal{G}^{op}}_{n_p}N_{n_p}$ for some map $\square[\underline{n}]\to K$ with $\underline{n}=(n_1,\dots,n_p) \in \seq(n)$. The same argument as above yields an element of $\P_{\alpha,\beta}^1[K]_{reg}$. This means that the left-hand map $\P_{\alpha,\beta}^1[K]_{reg}\to \P_{\alpha,\beta}^1||K||^{\mathcal{G}}$ is onto. It implies that the map $\P_{\alpha,\beta}^1[K]_{reg}\to \P_{\alpha,\beta}^1||K||^{\mathcal{G}}$ is a homeomorphism, thus the map $\P_{\alpha,\beta}^\ell[K]_{reg}\to \P_{\alpha,\beta}^\ell||K||^{\mathcal{G}}$ is a homeomorphism for all $\ell>0$. The proof is complete.
\epf

\section{Space of tame regular \mins{d}-paths and m-cofibrancy}
\label{main}

Some additional information about cube chains is required before proving Theorem~\ref{iso_reg_reg0}.

\bp \label{G-nonG}
Let $K$ be a precubical set. Let $(\alpha,\beta)\in K_0\p K_0$. Then there is the isomorphism of $\mathcal{G}$-spaces
\[
\mathbb{F}^{\mathcal{G}^{op}}_{n}\big(\liminj \D_{\alpha,\beta}(K,n)\big) \iso  \liminj \D_{\alpha,\beta}^{\mathcal{G}}(K,n)
\]
for all integers $n\geq 1$.
\ep

\bpf Let $\underline{n} = (n_1,\dots,n_p)$ and $\sum_i n_i=n$. Then 
\begin{align*}
	\D_{\alpha,\beta}^{\mathcal{G}}(K,n)(\square[\underline{n}]\to K) &= \mathbb{F}^{\mathcal{G}^{op}}_{n_1} N_{n_1} \ot \dots \ot \mathbb{F}^{\mathcal{G}^{op}}_{n_p} N_{n_p}\\
	&\iso \mathbb{F}^{\mathcal{G}^{op}}_{n}(N_{n_1}\p\dots \p N_{n_p}) \\
	&=\mathbb{F}^{\mathcal{G}^{op}}_{n}\D_{\alpha,\beta}(K,n)(\square[\underline{n}]\to K),
\end{align*}
the first equality by definition of $\D_{\alpha,\beta}^{\mathcal{G}}(K,n)$, the isomorphism by Proposition~\ref{Ftenseur} and the last equality by definition of $\D_{\alpha,\beta}^{\mathcal{G}}(K,n)$. The conclusion follows from Proposition~\ref{ev-adj}.
\epf

\bp \label{cp3}
For all precubical sets $K$, there is the natural isomorphism of flows \[\lmoore(||K||^{\mathcal{G}}) \iso ||K||.\]
\ep

\bpf
Let $(\alpha,\beta)\in K_0\p K_0$. There is the sequence of homeomorphisms 
\begin{align*}
	\P_{\alpha,\beta}\lmoore(||K||^{\mathcal{G}})
	&\iso \liminj \big(\P_{\alpha,\beta} ||K||^{\mathcal{G}}\big)\\
	& \iso \coprod_{n\geq 1}\liminj \bigg(\mathbb{F}^{\mathcal{G}^{op}}_{n}\big(\liminj \D_{\alpha,\beta}(K,n)\big)\bigg) \\
	& \iso \coprod_{n\geq 1} \liminj \D_{\alpha,\beta}(K,n)\\
	&\iso \P_{\alpha,\beta} ||K||,
\end{align*}
the first homeomorphism by definition of $\lmoore$, the second homeomorphism by Proposition~\ref{G-nonG} and since colimits commute with coproducts, the third homeomorphism by \cite[Proposition~5.8]{dgrtop}, and the last homeomorphism by definition of the flow $||K||$.
\epf

\begin{nota} \cite[Notation~A.1]{NaturalRealization}
	Let $n\geq 3$. Let $\mathcal{B}_n$ be the set of precubical sets $A$ such that $A\subset \de\square[n]$ and such that $|A|_{geom} \subset [0,1]^n$ contains a $d$-path of $\vec{N}^n_{0_n,1_n}(\square[n])$ which does not intersect $\{0,1\}^n\backslash\{0_n,1_n\}$. In particular, this means that $0_n,1_n$ are two vertices of $A$. One has $\de\square[n]\in \mathcal{B}_n$.
\end{nota}

\bd \label{carac_spatial}
A precubical set is \textit{spatial} if it is orthogonal to the set of maps of precubical sets 
\[\bigg\{\square[n]\sqcup_A \square[n]\longrightarrow \square[n]\mid n\geq 3 \hbox{ and }A\in \mathcal{B}_n\bigg\}.\]
\ed

Every proper precubical set in the sense of \cite[page~499]{MR3722069} is spatial by \cite[Proposition~7.5]{NaturalRealization}. In particular, for all $n\geq 0$, the precubical sets $\de\square[n]$ and $\square[n]$ are spatial, as well as all geometric precubical sets in the sense of \cite[Definition~1.18]{zbMATH07226006} and all non-positively curved precubical sets in the sense of \cite[Definition~1.28]{zbMATH07226006}, since they are proper. Also every $2$-dimensional precubical set is spatial by \cite[Corollary~A.3]{NaturalRealization}.

\bth \label{iso_reg_reg0} 
For all precubical sets $K$, there is a natural weak equivalence of the h-model structure of Moore flows \[[K]_{reg} \longrightarrow \moore^{\mathcal{G}}(|K|^t_{reg}).\] Moreover, the weak equivalence above is an isomorphism of Moore flows if and only if $K$ is spatial.
\eth

\bpf
From the cocone 
\[
\bigg(\moore^{\mathcal{G}}|\square[n]|^t_{reg}\bigg)_{\square[n]\to K} \stackrel{\bullet}\longrightarrow \moore^{\mathcal{G}}|K|^t_{reg}
\] 
we deduce the natural map of Moore flows $g:[K]_{reg} \longrightarrow \moore^{\mathcal{G}}|K|^t_{reg}$. It is bijective on states. Let $(\alpha,\beta)\in K_0\p K_0$. Consider the following commutative diagram of topological spaces
\[
\xymatrix@C=2.7em@R=3em
{
	\ar@{->}[d]_-{h_1}\displaystyle\coprod_{n\geq 1} \mathcal{G}(1,n) \p \liminj \D_{\alpha,\beta}(K,n) \ar@{->}[rr]^-{h_3} && \displaystyle\coprod_{n\geq 1} \mathcal{G}(1,n) \p \vec{N}^n_{\alpha,\beta}(K)\ar@{->}[d]^-{\Phi^1}\\
	\P_{\alpha,\beta}^1 ||K||^{\mathcal{G}} \fr{h_2}& \P_{\alpha,\beta}^1 [K]_{reg} \fr{\P_{\alpha,\beta}^1g}& \P_{\alpha,\beta}^1\moore^{\mathcal{G}}|K|^t_{reg} = \vec{R}^1_{\alpha,\beta}(K)
}
\]
where 1) the map $\Phi^1$ is the homeomorphism of Theorem~\ref{Psi}, 2) the map $h_1$ is the homeomorphism of Proposition~\ref{G-nonG}, 3) the map $h_2$ is the homeomorphism given by Theorem~\ref{doublebarG}, and 4) the map $h_3$ is induced by the continuous map $\D_{\alpha,\beta}(K,n)\to \vec{N}^n_{\alpha,\beta}(K)$ given by the Moore composition of tame natural $d$-paths defined as follows. Let $\underline{n}=(n_1,\dots,n_p)$, $n=\sum_i n_i$ and consider a map $\square[\underline{n}] \to K \in \Ch_{\alpha,\beta}(K,n)$. It gives rise to a sequence of cubes $(c_1,\dots,c_p)$ of $K$. The continuous map $\D_{\alpha,\beta}(K,n)\to \vec{N}^n_{\alpha,\beta}(K)$ takes \[(\gamma_1,\dots,\gamma_p)\in \D_{\alpha,\beta}(K,n)(\square[\underline{n}]\to K) = N_{n_1} \p \dots \p N_{n_p}\] to the Moore composition of tame natural $d$-paths \[(|c_1|_{geom}\gamma_1)*\dots *(|c_p|_{geom}\gamma_p)\in\vec{N}^n_{\alpha,\beta}(K)\] which is denoted by $[c_1;\gamma_1]*\dots [c_p;\gamma_p]$ in \cite{NaturalRealization}. By \cite[Theorem~7.7]{NaturalRealization}, the map $h_3$ is a homeomorphism when $K$ is spatial. For a general precubical set $K$, the map $h_3$ is a homotopy equivalence by \cite[Theorem~7.8]{NaturalRealization}. Hence the first part of the proof is complete.

Conversely, suppose that the natural map $[K]_{reg} \to \moore^{\mathcal{G}}(|K|^t_{reg})$ is an isomorphism of Moore flows. Then there is the sequence of isomorphisms of flows 
\[
||K||\iso \lmoore||K||^{\mathcal{G}}\iso \lmoore[K]_{reg} \stackrel{\iso} \longrightarrow \lmoore\moore^{\mathcal{G}}(|K|^t_{reg})\iso \Cat(|K|^t_{reg}) \iso |K|_{tc},
\]
the first isomorphism by Proposition~\ref{cp3}, the second isomorphism by Theorem~\ref{doublebarG}, the third isomorphism by hypothesis, the fourth isomorphism by Theorem~\ref{decomposing2} and the last isomorphism by Proposition~\ref{same}. This isomorphism of flows from $||K||$ to $|K|_{tc}$ is the identity of $K_0$ on states, and for all $(\alpha,\beta)\in K_0\p K_0$, it takes an element of $\P_{\alpha,\beta}||K|| \iso \coprod_{n\geq 1} \liminj \D_{\alpha,\beta}(K,n)$ to an element of $\P_{\alpha,\beta}|K|_{tc} = \coprod_{n\geq 1} \vec{N}^n_{\alpha,\beta}(K)$ as described above. By \cite[Theorem~A.2]{NaturalRealization}, the precubical set $K$ is spatial.
\epf

A by-product of the calculation made in the proof of Theorem~\ref{iso_reg_reg0} is the following fact. The image by the functor $\lmoore:\dtopG\to \dtop$ of the map $[K]_{reg} \longrightarrow \moore^{\mathcal{G}}(|K|^t_{reg})$ is the map of flows $||K||\to |K|_{tc}$ which is induced by the maps $N_{n_1} \p \dots \p N_{n_p} \to \vec{N}^{n}_{\alpha,\beta}(K)$ for $\square[\underline{n}] \to K \in \Ch_{\alpha,\beta}(K,n)$ with $\underline{n}=(n_1,\dots,n_p)$ and $n=\sum_i n_i$ as described in the core of the proof. It is a weak equivalence of the h-model structure of flows, and an isomorphism if and only if $K$ is spatial by \cite[Theorem~7.6 and Theorem~7.7]{NaturalRealization}.

\bp \label{pushout}
For all $n\geq 1$, there is the pushout diagram of Moore flows 
\[
\xymatrix@C=3em@R=3em
{
	\glob(\mathbb{F}^{\mathcal{G}^{op}}_n\de N_n) \fd{} \fr{} & [\de\square[n]]_{reg} \fd{} \\
	\glob(\mathbb{F}^{\mathcal{G}^{op}}_n N_n) \fr{} & \cocartesien [\square[n]]_{reg}
}
\]
\ep

\bpf
By Corollary~\ref{cp2} and $\square[n]$ and $\de\square[n]$ being spatial for all $n\geq 0$, we obtain using Theorem~\ref{iso_reg_reg0} the commutative diagram of $\mathcal{G}$-spaces 
\[
\xymatrix@C=3em@R=3em
{
	\mathbb{F}^{\mathcal{G}^{op}}_n\de N_n \fd{} \fr{\iso} & \P_{0_{n},1_{n}}[\de\square[n]]_{reg} \fd{} \\
	\mathbb{F}^{\mathcal{G}^{op}}_n N_n \fr{\iso} &  \P_{0_{n},1_{n}}[\square[n]]_{reg}
}
\]
where the two horizontal maps are isomorphisms of $\mathcal{G}$-spaces. This commutative diagram of $\mathcal{G}$-spaces is therefore a pushout diagram and the proof is complete.
\epf

\bth \label{m-cof-moore-flow}
For all precubical sets $K$, the Moore flow $[K]_{reg}$ is m-cofibrant.
\eth

\bpf A map of Moore flows $f:X\to Y$ satisfies the RLP with respect to the map of Moore flows $\glob(\mathbb{F}^{\mathcal{G}^{op}}_n\de N_n) \to \glob(\mathbb{F}^{\mathcal{G}^{op}}_n N_n)$ if and only if for each $(\alpha,\beta)\in K_0\p K_0$, the map of $\mathcal{G}$-spaces $f:\P_{\alpha,\beta}X \to \P_{f(\alpha),f(\beta)}Y$ satisfies the RLP with respect to the map of $\mathcal{G}$-spaces $\mathbb{F}_n^{\mathcal{G}^{op}}\de N_n \to \mathbb{F}_n^{\mathcal{G}^{op}} N_n$. Therefore, by Proposition~\ref{ev-adj}, a map of Moore flows $f:X\to Y$ satisfies the RLP with respect to the map of Moore flows $\glob(\mathbb{F}^{\mathcal{G}^{op}}_n\de N_n) \to \glob(\mathbb{F}^{\mathcal{G}^{op}}_n N_n)$ if and only the map of topological spaces $\P^n_{\alpha,\beta}X \to \P^n_{f(\alpha),f(\beta)}Y$ satisfies the RLP with respect to the map of topological spaces $\de N_n \to N_n$. By \cite[Theorem~5.9]{NaturalRealization}, the latter map is an m-cofibration of spaces. By Theorem~\ref{qhmMooreFlow}, the trivial fibrations of the m-model structures of Moore flows are objectwise. We then deduce that the map of Moore flows $\glob(\mathbb{F}^{\mathcal{G}^{op}}_n\de N_n) \to \glob(\mathbb{F}^{\mathcal{G}^{op}}_n N_n)$ is an m-cofibration of Moore flows for all $n\geq 1$. Using Proposition~\ref{pushout}, and since $[\square[0]]_{reg} = \{0\}$ (the Moore flow without execution paths and one state $0$) is m-cofibrant, we deduce that for all precubical sets $K$, the Moore flow $[K]_{reg}$ is m-cofibrant.  
\epf

\begin{cor} \label{space-mcof}
	Let $K$ be a precubical set. Let $(\alpha,\beta)\in K_0\p K_0$. The space of tame regular $d$-paths from $\alpha$ to $\beta$ in the geometric realization of $K$ is homotopy equivalent to a CW-complex.
\end{cor}

\bpf
The Moore flow $\moore^{\mathcal{G}}[K]_{reg}$ is m-cofibrant by Theorem~\ref{m-cof-moore-flow}. Thus it is weakly equivalent in the h-model structure of Moore flows to a q-cofibrant Moore flow $X$ by \cite[Corollary~3.7]{mixed-cole}. The latter has projective q-cofibrant $\mathcal{G}$-spaces of execution paths by \cite[Theorem~9.11]{Moore1}. It implies that $\P_{\alpha,\beta}X$ is injective m-cofibrant by \cite[Corollary~7.2]{dgrtop}. Thus $\P^1_{\alpha,\beta}X$ is an m-cofibrant topological space. Since weak equivalences of the projective h-model structure of $\mathcal{G}$-spaces are objectwise homotopy equivalences, it implies that $\P^1_{\alpha,\beta}\moore^{\mathcal{G}}[K]_{reg}$ is homotopy equivalent to $\P^1_{\alpha,\beta}X$. By Theorem~\ref{iso_reg_reg0}, $\P^1_{\alpha,\beta}\moore^{\mathcal{G}}[K]_{reg}$ is homotopy equivalent to $\P^{top}_{\alpha,\beta}|K|^t_{reg}$. When $\alpha=\beta$, the space of tame regular $d$-paths from $\alpha$ to $\beta$ is homeomorphic to the disjoint union of $\{\alpha\}$ and the space of nonconstant tame regular $d$-paths from $\alpha$ to $\beta$. For $\alpha\neq \beta$, the latter remark is pointless. Hence the proof is complete. 
\epf

Corollary~\ref{space-mcof} was previously established by several other authors by explicitly constructing homotopy equivalences with CW-complexes. At first, the adjective regular can be removed from the statement of Corollary~\ref{space-mcof} thanks to a result due to Raussen \cite[Proposition~2.16]{MR2521708}: see the proof of Theorem~\ref{comparison-tame-nottame-Moore-flow}. The first method to establish Corollary~\ref{space-mcof} without the adjective regular is due to Ziemia\'{n}ski \cite[Theorem~6.1, Theorem~7.5 and Theorem~7.6]{MR4070250}. It consists of proving that the space of tame $d$-paths from $\alpha$ to $\beta$ in the geometric realization is homotopy equivalent to the classifying space of a small category, namely the small category of cube chains from $\alpha$ to $\beta$ (as already pointed before \cite[Corollary~7.9]{NaturalRealization}, the word ``weak homotopy equivalence'' can be replaced by ``homotopy equivalence'' in \cite[Theorem~7.5 and Theorem~7.6]{MR4070250}). Another method is proposed by Paliga and Ziemia\'{n}ski in \cite[Proposition~6.5]{PZ}, which uses another result due to Raussen, namely \cite[Proposition~3.15]{MR2521708}. However, it works only for finite precubical sets. Finally, another method is expounded by Raussen in \cite{raussen2021strictifying} which works only for proper non-self-linked precubical sets.

The reader may wonder whether the multipointed $d$-space $|K|^t_{reg}$ is m-cofibrant as well for all precubical sets $K$. We are not even able to prove that it is always h-cofibrant. We can only prove what follows:

\bp \label{m-spatial}
Let $K$ be a spatial precubical set. Then the multipointed $d$-space $|K|^t_{reg}$ is m-cofibrant if it is h-cofibrant.
\ep

\bpf
The functor $\moore^{\mathcal{G}}:\ptop{\mathcal{G}}\to \dtopG$ takes (trivial resp.) m-fibrations of multipointed $d$-spaces to (trivial resp.) m-fibrations of Moore flows by definition of the m-model structures (see Theorem~\ref{three} and Theorem~\ref{qhmMooreFlow}). Thus, the functor $\moore^{\mathcal{G}}:\ptop{\mathcal{G}}\to \dtopG$ is a right Quillen adjoint between the m-model structures. We have the commutative diagram of right Quillen adjoints 
\[
\xymatrix@C=3em@R=3em
{
	\ptop{\mathcal{G}}_m \ar@{->}[d]_-{\id_{\ptop{\mathcal{G}}}}\fr{\moore^{\mathcal{G}}} & \dtopG_m \ar@{->}[d]^-{\id_{\dtopG}}\\
	\ptop{\mathcal{G}}_q \ar@{-}[r]_-{\moore^{\mathcal{G}}} & \dtopG_q
}
\]
The bottom horizontal arrow is a right Quillen equivalence by \cite[Theorem~8.1]{Moore2}. It is a general fact about mixed model structures that the two vertical arrows are right Quillen equivalences. Thus, the top right Quillen adjoint is a right Quillen equivalence. All multipointed $d$-spaces are m-fibrant by Theorem~\ref{three}. Since the Moore flow $[K]_{reg}$ is m-cofibrant by Theorem~\ref{m-cof-moore-flow}, we deduce that the counit map 
\[
\lmoore^{\mathcal{G}}(\moore^{\mathcal{G}}|K|^t_{reg}) \iso  \lmoore^{\mathcal{G}}([K]_{reg})\longrightarrow |K|^t_{reg},
\]
is a weak equivalence of the m-model structure of multipointed $d$-spaces by Theorem~\ref{iso_reg_reg0}, the precubical set $K$ being spatial by hypothesis. Moreover, $\lmoore^{\mathcal{G}}$ being a left Quillen adjoint, the multipointed $d$-space $\lmoore^{\mathcal{G}}(\moore^{\mathcal{G}}|K|^t_{reg})$ is m-cofibrant. For all $(\alpha,\beta)\in K_0\p K_0$, we deduce the weak homotopy equivalence 
\[
\P^{top}_{\alpha,\beta}\lmoore^{\mathcal{G}}(\moore^{\mathcal{G}}|K|^t_{reg}) \longrightarrow \P^{top}_{\alpha,\beta}|K|^t_{reg}.
\]
By \cite[Theorem~8.6]{QHMmodel}, the topological space $\P^{top}_{\alpha,\beta}\lmoore^{\mathcal{G}}(\moore^{\mathcal{G}}|K|^t_{reg})$ is m-cofibrant. By Corollary~\ref{space-mcof}, the space $\P^{top}_{\alpha,\beta}|K|^t_{reg}$ is an m-cofibrant space. Thus the weak homotopy equivalence 
\[
\P^{top}_{\alpha,\beta}\lmoore^{\mathcal{G}}(\moore^{\mathcal{G}}|K|^t_{reg}) \longrightarrow \P^{top}_{\alpha,\beta}|K|^t_{reg}.
\]
is a weak homotopy equivalence between m-cofibrant spaces. By \cite[Corollary~3.4]{mixed-cole}, the latter map is therefore a homotopy equivalence. This means that the counit map $\lmoore^{\mathcal{G}}(\moore^{\mathcal{G}}|K|^t_{reg}) \to |K|^t_{reg}$ is a weak equivalence of the h-model structure of $\ptop{\mathcal{G}}$ when $K$ is spatial. This implies by \cite[Corollary~3.7]{mixed-cole} that the multipointed $d$-space $|K|^t_{reg}$ is m-cofibrant if it is h-cofibrant, the multipointed $d$-space $\lmoore^{\mathcal{G}}(\moore^{\mathcal{G}}|K|^t_{reg})$ being m-cofibrant.
\epf

The example of non h-cofibrant multipointed $d$-space provided in \cite[Proposition~6.19]{QHMmodel} suggests that the h-cofibrant objects are the objects without algebraic relations. According to this intuition, we have the following conjecture:

\begin{conj}
	For all precubical sets $K$, the multipointed $d$-space $|K|^t_{reg}$ is h-cofibrant. 
\end{conj}

\section{The regular realization of a precubical set}
\label{question-r}

The \textit{regular realization} of a precubical set is defined as follows: 

\bd \label{reg_rea}
Let $K$ be a precubical set. The \textit{regular realization} $|K|_{reg}$ of $K$ is the multipointed $d$-space having  the underlying space $|K|_{geom}$, the set of states $X^0$ and such that the set of execution paths from $\alpha$ to $\beta$ consists of the nonconstant regular $d$-paths from $\alpha$ to $\beta$ in the geometric realization of $K$.
\ed

The multipointed $d$-space $|K|_{reg}$ is well defined for all precubical sets $K$ because the normalized composition of two regular $d$-paths in the geometric realization of a precubical set is regular and because the reparametrization of a regular $d$-path by a homeomorphism is regular. There is an isomorphism of cocubical multipointed $d$-spaces $|\square[*]|^t_{reg}\iso |\square[*]|_{reg}$ by Proposition~\ref{cocubical-reg-cube}. By the universal property of the colimits, we obtain a natural map of multipointed $d$-spaces \[|K|^t_{reg}\longrightarrow |K|_{reg}.\] 

\bth \label{comparison-tame-nottame-Moore-flow}
For all precubical sets $K$, there is a natural weak equivalence for the h-model structure of Moore flows 
\[
\moore^{\mathcal{G}}(|K|^t_{reg}) \longrightarrow \moore^{\mathcal{G}}(|K|_{reg}).
\]
\eth

\bpf
The image by the functor $\moore^{\mathcal{G}}$ of the natural map of multipointed $d$-spaces $|K|^t_{reg} \to |K|_{reg}$ preserves the set of states. For each $(\alpha,\beta)\in K_0$, the continuous map \[\P_{\alpha,\beta}^1\moore^{\mathcal{G}}|K|^t_{reg} \subset \P_{\alpha,\beta}^1\moore^{\mathcal{G}}|K|_{reg}\] is the inclusion of the space of nonconstant tame regular $d$-paths from $\alpha$ to $\beta$ into the space of nonconstant regular $d$-paths from $\alpha$ to $\beta$, both equipped with the $\Delta$-kellefication of the compact-open topology. By \cite[Proposition~2.16]{MR2521708}, the inclusion of the space of regular $d$-paths from $\alpha$ to $\beta$ into the space of $d$-paths from $\alpha$ to $\beta$ has an inverse up to homotopy. The latter is a reparametrization of Moore $d$-paths (it is the composite of two reparametrizations, the normalization and the naturalization: see Raussen's proof). It therefore preserves tameness. Moreover, the homotopies with the identities are homotopies between reparametrizations: see Part~3 of the proof of \cite[Proposition~2.16]{MR2521708}. This implies that they preserve tameness as well. This implies that the inclusion of the space of tame regular $d$-paths from $\alpha$ to $\beta$ into the space of tame $d$-paths from $\alpha$ to $\beta$ has an inverse up to homotopy as well. The paper \cite{MR2521708} works with the compact-open topology, which is not an issue thanks to Proposition~\ref{homotopy-deltakelleyfication}. Besides, working with nonconstant $d$-paths only is not an issue either in the case of precubical sets: indeed either $\alpha\neq \beta$ and there is no constant $d$-paths from $\alpha$ to $\beta$ or $\alpha=\beta$ and the unique constant $d$-path is in a distinct path-connected component. Thus, the continuous map $\P_{\alpha,\beta}^1\moore^{\mathcal{G}}|K|^t_{reg} \subset \P_{\alpha,\beta}^1\moore^{\mathcal{G}}|K|_{reg}$ is homotopic to the inclusion of the space of nonconstant tame $d$-paths from $\alpha$ to $\beta$ into the space of nonconstant $d$-paths from $\alpha$ to $\beta$: the adjective regular can be omitted. Thanks to the tamification theorem \cite[Theorem~6.1]{MR4070250} and the \ttt\  for homotopy equivalences, we then obtain that the continuous map $\P_{\alpha,\beta}^1\moore^{\mathcal{G}}|K|^t_{reg} \subset \P_{\alpha,\beta}^1\moore^{\mathcal{G}}|K|_{reg}$ is a homotopy equivalence. We obtain a weak equivalence of the projective h-model structure of $\mathcal{G}$-spaces \[\P_{\alpha,\beta}\moore^{\mathcal{G}}|K|^t_{reg} \subset \P_{\alpha,\beta}\moore^{\mathcal{G}}|K|_{reg}\] for each $(\alpha,\beta)\in K_0\p K_0$. Hence the proof is complete.
\epf

\begin{figure}
	\begin{tikzpicture}
			\draw (0,0) -- (0,3) -- (3,3) -- (3,0) -- (0,0);
			\draw (0,3) -- (-1,4) -- (-1,1) -- (0,0);
			\draw (-1,4) -- (2,4) -- (3,3);
			\draw[dashed] (-1,1) -- (2,1) -- (3,0); 
			\draw[dashed] (2,1) -- (2,4); 
			\draw[->, line width=0.4mm, color=dark-red] (0,0) -- (1.8,3);
			\draw[->, line width=0.4mm, color=dark-red] (1.8,3) -- (2,4);
		\end{tikzpicture}
	\caption{Non-tame $d$-path from $0_3$ to $1_3$ in the boundary of the $3$-cube}
	\label{nontame}
\end{figure}
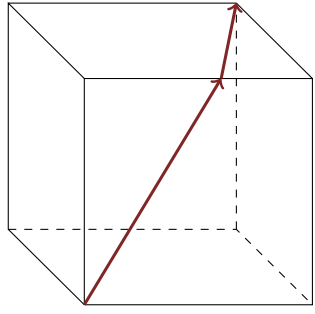

There is no reason for the weak equivalence $\moore^{\mathcal{G}}(|K|^t_{reg}) \to \moore^{\mathcal{G}}(|K|_{reg})$ to become an isomorphism when $K$ is spatial: the precubical set $\de\square[3]$ is spatial by \cite[Corollary~A.3]{NaturalRealization} and Figure~\ref{nontame} proves that the map $\moore^{\mathcal{G}}(|\de\square[3]|^t_{reg}) \to \moore^{\mathcal{G}}(|\de\square[3]|_{reg})$ is not an isomorphism.

\begin{cor} \label{comparison-tame-nottame-Moore-flow-dspace}
For all precubical sets $K$, there is a natural weak equivalence for the h-model structure of multipointed $d$-spaces 
\[
|K|^t_{reg} \longrightarrow |K|_{reg}.
\]
\end{cor}

\bpf
A straightforward consequence of the definition of a weak equivalence in the h-model structures of multipointed $d$-spaces and Moore flows is that the functor $\moore^{\mathcal{G}}:\ptop{\mathcal{G}}$ to $\dtopG$ reflects weak equivalences of the h-model structures (it is also true for the q-model structures and for the m-model structures). Thanks to Theorem~\ref{comparison-tame-nottame-Moore-flow}, the proof is complete.
\epf

\begin{cor} \label{reg_weak}
	For all precubical sets $K$, there exists a natural weak equivalence of the h-model structure of Moore flows 
	\[
	[K]_{reg} \longrightarrow \moore^{\mathcal{G}}(|K|_{reg}).
	\]
\end{cor}

\bpf
It is a consequence of Theorem~\ref{iso_reg_reg0} and Theorem~\ref{comparison-tame-nottame-Moore-flow}.
\epf

Corollary~\ref{space-mcof} gives a pure model categorical proof of the fact that the space of \textit{tame} regular $d$-paths between two vertices of a precubical set in its geometric realization is homotopy equivalent to a CW-complex. On the other hand, we do not know any pure model categorical proof of the same fact for the space of \textit{all} regular $d$-paths between two vertices. We can only deduce it from Corollary~\ref{space-mcof} and Theorem~\ref{comparison-tame-nottame-Moore-flow}, the latter relying on results from Raussen and Ziemia{\'{n}}ski. 

The following observations enable us to better understand the interaction between model categories and tameness. There exists a combinatorial model structure on multipointed $d$-spaces such that the multipointed $d$-spaces $|K|^t_{reg}$ are cofibrant for all precubical sets $K$. Indeed, it suffices to consider the minimal model structure in the sense of \cite[Theorem~1.4]{henry2020minimal} with respect to the set of maps $\{|\de\square[n]|^t_{reg}\subset |\square[n]|^t_{reg}\mid n\geq 0\}$ with or without the additional map $R:\{0,1\}\to \{0\}$. Note that the multipointed $d$-spaces $|K|^t_{reg}$ are in general not q-cofibrant for the q-model structure. In the same way, we could consider the minimal model structure on multipointed $d$-spaces with respect to the set $\mathcal{I}$ of maps $\{|\de\square[n]|_{reg}\subset |\square[n]|_{reg}\mid n\geq 0\}$ with or without the additional map $R:\{0,1\}\to \{0\}$. However, it turns out that for all $n\geq 3$, $|\de\square[n]|_{reg}$ is not $\mathcal{I}$-cofibrant by Proposition~\ref{nontractable}. This implies that \cite[Theorem~1.4]{henry2020minimal} cannot even be used to prove its existence. The best that can be said is that the latter model structure exists by assuming Vop\v{e}nka's principle thanks to \cite[Theorem~2.2]{rotho}. 

Before proving Proposition~\ref{nontractable}, let us start by a remark about $\mathcal{I}$-cellular multipointed $d$-spaces. The underlying space of an $\mathcal{I}$-cellular multipointed $d$-space is not necessarily the geometric realization of a precubical set because the attaching maps are not necessarily images of maps of precubical sets. However, the notion of tameness can be easily generalized to them. Consider a pushout diagram of multipointed $d$-spaces of the form 
\[
\xymatrix@C=3em@R=3em
{
	|\de\square[n]|_{reg} \fd{}\fr{} & X \fd{}\\
	|\square[n]|_{reg} \fr{} & \cocartesien Y
}
\]
with $n\geq 0$ such that $X$ is $\mathcal{I}$-cellular and such that all its execution paths are tame. By definition of $|\square[n]|_{reg} = |\square[n]|^t_{reg}$ for all $n\geq 0$ (see Proposition~\ref{cocubical-reg-cube}), and since colimits are calculated by taking a final structure as explained at the end of Section~\ref{reg_sec}, we deduce that all execution paths of $Y$ are tame even if $|\de\square[n]|_{reg}$ may contain non-tame execution paths: for example the non-tame $d$-path depicted in Figure~\ref{nontame} becomes tame if $\de\square[3]$ is embedded in $\square[3]$. By an immediate transfinite induction, we obtain that all execution paths of an $\mathcal{I}$-cellular multipointed $d$-space are tame. The converse is false. Indeed, the multipointed $d$-space $|\square[3]|_{reg}$ is not $\mathcal{I}$-cellular, $|\de\square[3]|_{reg}$ being not $\mathcal{I}$-cellular, and all its execution paths are tame. One can also easily prove by transfinite induction on the $\mathcal{I}$-cellular decomposition of an $\mathcal{I}$-cellular multipointed $d$-space $X$ that $X^0$ is a discrete subspace of $|X|$ and that $|X|$ is Hausdorff. Recall that, for a general multipointed $d$-space $X$, $X^0$ is not necessarily discrete: see Figure~\ref{contracting}. 

\bp \label{nontractable}
Let $n\geq 3$. The multipointed $d$-space $|\de\square[n]|_{reg}$ is not $\mathcal{I}$-cofibrant.
\ep

\bpf[Sketch of proof]
The multipointed $d$-space $|\de\square[n]|_{reg}$ is not $\mathcal{I}$-cellular precisely because it contains non-tame execution paths (see Figure~\ref{nontame} for the case $n=3$). Suppose that $|\de\square[n]|_{reg}$ is a retract of an $\mathcal{I}$-cellular multipointed $d$-space $X$. This means that the identity of $|\de\square[n]|_{reg}$ factors as a composite \[|\de\square[n]|_{reg} \longrightarrow X \longrightarrow |\de\square[n]|_{reg}.\] The left-hand map being one-to-one and $|X|$ being Hausdorff, this implies that the compact $|\de\square[n]|_{geom}$, which is homeomorphic to the $(n-1)$-dimensional sphere $\mathbf{S}^{n-1}$, is homeomorphic to its compact image $K$ in $|X|$. The composite map of groups \[\mathbb{Z}\iso \pi_{n-1}(|\de\square[n]|_{geom}) \longrightarrow \pi_{n-1}(|X|) \longrightarrow \pi_{n-1}(|\de\square[n]|_{geom}) \iso \mathbb{Z}\] being the identity of $\mathbb{Z}$ ($\pi_{n-1}$ denotes the $(n-1)$-th homotopy group and we choose any point of $|\de\square[n]|_{geom}$ as a base point), we deduce that the compact $K\subset |X|$ is not contained in a $p$-cube for some $p\geq n$: intuitively the hole cannot be filled in $|X|$. Now, recall that every execution path of $|\de\square[n]|_{reg}$ is taken to an execution path of $X$ which is necessarily tame because, by hypothesis, $X$ is $\mathcal{I}$-cellular. By deforming continuously a non-tame $d$-path from $0_n$ to $1_n$ in $|\de\square[n]|_{geom}$ to obtain a family $(\gamma_u)_{u\in [0,1]}$ of non-tame $d$-paths from $0_n$ to $1_n$ in $|\de\square[n]|_{geom}$ such that $u\neq u'$ implies $\gamma_u([0,1]) \cap \gamma_{u'}([0,1]) =\{0_n,1_n\}$, we then build an infinite set of elements of $X^0$ in $K$, the left-hand map $|\de\square[n]|_{geom} \to |X|$ being one-to-one. However, $X^0 \cap K$ is finite, $X^0$ being discrete. It is a contradiction. 
\epf

%

\begin{thebibliography}{10}
	
	\bibitem{topologicalcat}
	J.~Ad{\'a}mek, H.~Herrlich, and G.~E. Strecker.
	\newblock Abstract and concrete categories: the joy of cats.
	\newblock {\em Repr. Theory Appl. Categ.}, (17):1--507, 2006.
	\newblock Reprint of the 1990 original [Wiley, New York; MR1051419].
	
	\bibitem{TheBook}
	J.~Ad{\'a}mek and J.~Rosick{\'y}.
	\newblock {\em Locally presentable and accessible categories}.
	\newblock Cambridge University Press, Cambridge, 1994.
	\newblock \href {https://doi.org/10.1017/cbo9780511600579.004}
	{\path{https://doi.org/10.1017/cbo9780511600579.004}}.
	
	\bibitem{Barthel-Riel}
	T.~Barthel and E.~Riehl.
	\newblock On the construction of functorial factorizations for model
	categories.
	\newblock {\em Algebr. Geom. Topol.}, 13(2):1089--1124, 2013.
	\newblock \href {https://doi.org/10.2140/agt.2013.13.1089}
	{\path{https://doi.org/10.2140/agt.2013.13.1089}}.
	
	\bibitem{Borceux2}
	F.~Borceux.
	\newblock {\em Handbook of categorical algebra. 2}.
	\newblock Cambridge University Press, Cambridge, 1994.
	\newblock Categories and structures.
	\newblock \href {https://doi.org/10.1017/cbo9780511525865}
	{\path{https://doi.org/10.1017/cbo9780511525865}}.
	
	\bibitem{Brown_cube}
	R.~Brown and P.~J. Higgins.
	\newblock On the algebra of cubes.
	\newblock {\em J. Pure Appl. Algebra}, 21(3):233--260, 1981.
	\newblock \href {https://doi.org/10.1016/0022-4049(81)90018-9}
	{\path{https://doi.org/10.1016/0022-4049(81)90018-9}}.
	
	\bibitem{mixed-cole}
	M.~Cole.
	\newblock Mixing model structures.
	\newblock {\em Topology Appl.}, 153(7):1016--1032, 2006.
	\newblock \href {https://doi.org/10.1016/j.topol.2005.02.004}
	{\path{https://doi.org/10.1016/j.topol.2005.02.004}}.
	
	\bibitem{HomotopicalCategory}
	W.~G. Dwyer, P.~S. Hirschhorn, D.~M. Kan, and J.~H. Smith.
	\newblock {\em Homotopy limit functors on model categories and homotopical
		categories}, volume 113 of {\em Mathematical Surveys and Monographs}.
	\newblock American Mathematical Society, Providence, RI, 2004.
	\newblock \href {https://doi.org/10.1090/surv/113}
	{\path{https://doi.org/10.1090/surv/113}}.
	
	\bibitem{reparam}
	U.~Fahrenberg and M.~Raussen.
	\newblock Reparametrizations of continuous paths.
	\newblock {\em J. Homotopy Relat. Struct.}, 2(2):93--117, 2007.
	
	\bibitem{DAT_book}
	L.~Fajstrup, E.~Goubault, E.~Haucourt, S.~Mimram, and M.~Raussen.
	\newblock {\em Directed algebraic topology and concurrency. {With} a foreword
		by {Maurice} {Herlihy} and a preface by {Samuel} {Mimram}}.
	\newblock SpringerBriefs Appl. Sci. Technol. Springer, 2016.
	\newblock \href {https://doi.org/10.1007/978-3-319-15398-8}
	{\path{https://doi.org/10.1007/978-3-319-15398-8}}.
	
	\bibitem{FR}
	L.~Fajstrup and J.~Rosick{\'y}.
	\newblock A convenient category for directed homotopy.
	\newblock {\em Theory Appl. Categ.}, 21:7--20, 2008.
	
	\bibitem{MR1074175}
	R.~Fritsch and R.~A. Piccinini.
	\newblock {\em Cellular structures in topology}, volume~19 of {\em Cambridge
		Studies in Advanced Mathematics}.
	\newblock Cambridge University Press, Cambridge, 1990.
	\newblock \href {https://doi.org/10.1017/cbo9780511983948}
	{\path{https://doi.org/10.1017/cbo9780511983948}}.
	
	\bibitem{GKR18}
	R.~Garner, M.~{K\c{e}dziorek}, and E.~Riehl.
	\newblock {Lifting accessible model structures}.
	\newblock {\em {J. Topol.}}, 13(1):59--76, 2020.
	\newblock \href {https://doi.org/10.1112/topo.12123}
	{\path{https://doi.org/10.1112/topo.12123}}.
	
	\bibitem{model3}
	P.~Gaucher.
	\newblock A model category for the homotopy theory of concurrency.
	\newblock {\em Homology Homotopy Appl.}, 5(1):p.549--599, 2003.
	\newblock \href {https://doi.org/10.4310/hha.2003.v5.n1.a20}
	{\path{https://doi.org/10.4310/hha.2003.v5.n1.a20}}.
	
	\bibitem{realization}
	P.~{Gaucher}.
	\newblock {Globular realization and cubical underlying homotopy type of time
		flow of process algebra}.
	\newblock {\em {New York J. Math.}}, 14:101--137, 2008.
	
	\bibitem{mdtop}
	P.~Gaucher.
	\newblock Homotopical interpretation of globular complex by multipointed
	d-space.
	\newblock {\em Theory Appl. Categ.}, 22(22):588--621, 2009.
	
	\bibitem{dgrtop}
	P.~Gaucher.
	\newblock Enriched diagrams of topological spaces over locally contractible
	enriched categories.
	\newblock {\em New York J. Math.}, 25:1485--1510, 2019.
	
	\bibitem{Moore1}
	P.~Gaucher.
	\newblock Homotopy theory of {M}oore flows ({I}).
	\newblock {\em {Compositionality}}, 3(3), 2021.
	\newblock \href {https://doi.org/10.32408/compositionality-3-3}
	{\path{https://doi.org/10.32408/compositionality-3-3}}.
	
	\bibitem{Moore2}
	P.~Gaucher.
	\newblock Homotopy theory of {M}oore flows ({II}).
	\newblock {\em {Extr. Math.}}, 36(2):157--239, 2021.
	\newblock \href {https://doi.org/10.17398/2605-5686.36.2.157}
	{\path{https://doi.org/10.17398/2605-5686.36.2.157}}.
	
	\bibitem{leftproperflow}
	P.~Gaucher.
	\newblock Left properness of flows.
	\newblock {\em Theory Appl. Categ.}, 37(19):562--612, 2021.
	
	\bibitem{QHMmodel}
	P.~Gaucher.
	\newblock Six model categories for directed homotopy.
	\newblock {\em Categ. Gen. Algebr. Struct. Appl.}, 15(1):145--181, 2021.
	\newblock \href {https://doi.org/10.52547/cgasa.15.1.145}
	{\path{https://doi.org/10.52547/cgasa.15.1.145}}.
	
	\bibitem{NaturalRealization}
	P.~Gaucher.
	\newblock Comparing cubical and globular directed paths.
	\newblock {\em {Fund. Math.}}, 262(3):259--286, 2023.
	\newblock \href {https://doi.org/10.4064/fm219-3-2023}
	{\path{https://doi.org/10.4064/fm219-3-2023}}.
	
	\bibitem{Nonunital}
	P.~Gaucher.
	\newblock Comparing the non-unital and unital settings for directed homotopy.
	\newblock {\em Cah. Topol. G\'eom. Diff\'er. Cat\'eg.}, LXIV-2:176--197, 2023.
	
	\bibitem{Moore3}
	P.~Gaucher.
	\newblock Homotopy theory of {M}oore flows ({III}), 2023.
	\newblock \href {https://doi.org/10.48550/arXiv.2303.16174}
	{\path{https://doi.org/10.48550/arXiv.2303.16174}}.
	
	\bibitem{zbMATH07226006}
	E.~{Goubault} and S.~{Mimram}.
	\newblock {Directed homotopy in non-positively curved spaces}.
	\newblock {\em {Log. Methods Comput. Sci.}}, 16(3):55, 2020.
	\newblock Id/No 4.
	\newblock \href {https://doi.org/10.23638/LMCS-16(3:4)2020}
	{\path{https://doi.org/10.23638/LMCS-16(3:4)2020}}.
	
	\bibitem{henry2020minimal}
	S.~Henry.
	\newblock Minimal model structures, 2020.
	\newblock \href {https://doi.org/10.48550/arXiv.2011.13408}
	{\path{https://doi.org/10.48550/arXiv.2011.13408}}.
	
	\bibitem{HKRS17}
	K.~Hess, M.~K\c{e}dziorek, E.~Riehl, and B.~Shipley.
	\newblock A necessary and sufficient condition for induced model structures.
	\newblock {\em J. Topol.}, 10(2):324--369, 2017.
	\newblock \href {https://doi.org/10.1112/topo.12011}
	{\path{https://doi.org/10.1112/topo.12011}}.
	
	\bibitem{ref_model2}
	P.~S. Hirschhorn.
	\newblock {\em Model categories and their localizations}, volume~99 of {\em
		Mathematical Surveys and Monographs}.
	\newblock American Mathematical Society, Providence, RI, 2003.
	\newblock \href {https://doi.org/10.1090/surv/099}
	{\path{https://doi.org/10.1090/surv/099}}.
	
	\bibitem{MR99h:55031}
	M.~Hovey.
	\newblock {\em Model categories}.
	\newblock American Mathematical Society, Providence, RI, 1999.
	\newblock \href {https://doi.org/10.1090/surv/063}
	{\path{https://doi.org/10.1090/surv/063}}.
	
	\bibitem{KellyEnriched}
	G.~M. Kelly.
	\newblock Basic concepts of enriched category theory.
	\newblock {\em Repr. Theory Appl. Categ.}, 10:vi+137 pp., 2005.
	\newblock Reprint of the 1982 original [Cambridge Univ. Press, Cambridge].
	
	\bibitem{LawvereMetric}
	F.~W. Lawvere.
	\newblock Metric spaces, generalized logic, and closed categories.
	\newblock {\em Repr. Theory Appl. Categ}, 2002(1):1--37, 2002.
	
	\bibitem{MR1712872}
	S.~Mac~Lane.
	\newblock {\em Categories for the working mathematician}.
	\newblock Springer-Verlag, New York, second edition, 1998.
	\newblock \href {https://doi.org/10.1007/978-1-4757-4721-8}
	{\path{https://doi.org/10.1007/978-1-4757-4721-8}}.
	
	\bibitem{MR1300636}
	S.~Mac~Lane and I.~Moerdijk.
	\newblock {\em Sheaves in geometry and logic}.
	\newblock Universitext. Springer-Verlag, New York, 1994.
	\newblock A first introduction to topos theory, Corrected reprint of the 1992
	edition.
	\newblock \href {https://doi.org/10.1007/978-1-4612-0927-0}
	{\path{https://doi.org/10.1007/978-1-4612-0927-0}}.
	
	\bibitem{MoserLyne}
	L.~Moser.
	\newblock Injective and projective model structures on enriched diagram
	categories.
	\newblock 21(2):279--300, 2019.
	\newblock \href {https://doi.org/10.4310/hha.2019.v21.n2.a15}
	{\path{https://doi.org/10.4310/hha.2019.v21.n2.a15}}.
	
	\bibitem{PZ}
	J.~Paliga and K.~Ziemia\'nski.
	\newblock Configuration spaces and directed paths on the final precubical set.
	\newblock {\em {Fund. Math.}}, 257(3):229--263, 2022.
	\newblock \href {https://doi.org/10.4064/fm114-9-2021}
	{\path{https://doi.org/10.4064/fm114-9-2021}}.
	
	\bibitem{reparam-fixed}
	M.~Raussen.
	\newblock Reparametrizations with given stop data.
	\newblock {\em J. Homotopy Relat. Struct.}, 4(1):1--5, 2009.
	
	\bibitem{MR2521708}
	M.~Raussen.
	\newblock Trace spaces in a pre-cubical complex.
	\newblock {\em Topology Appl.}, 156(9):1718--1728, 2009.
	\newblock \href {https://doi.org/10.1016/j.topol.2009.02.003}
	{\path{https://doi.org/10.1016/j.topol.2009.02.003}}.
	
	\bibitem{Raussen2012}
	M.~Raussen.
	\newblock Simplicial models for trace spaces {II}: General higher dimensional
	automata.
	\newblock {\em Algebr. Geom. Topol.}, 12(3):1741--1761, August 2012.
	\newblock \href {https://doi.org/10.2140/agt.2012.12.1741}
	{\path{https://doi.org/10.2140/agt.2012.12.1741}}.
	
	\bibitem{raussen2021strictifying}
	M.~Raussen.
	\newblock Strictifying and taming directed paths in higher dimensional
	automata.
	\newblock {\em Math. Struct. Comput. Sci.}, 31(2):193--213, 2021.
	\newblock \href {https://doi.org/10.1017/s0960129521000128}
	{\path{https://doi.org/10.1017/s0960129521000128}}.
	
	\bibitem{MR2506258}
	J.~Rosick{\'y}.
	\newblock On combinatorial model categories.
	\newblock {\em Appl. Categ. Structures}, 17(3):303--316, 2009.
	\newblock \href {https://doi.org/10.1007/s10485-008-9171-2}
	{\path{https://doi.org/10.1007/s10485-008-9171-2}}.
	
	\bibitem{zbMATH06722019}
	J.~Rosick{\'y}.
	\newblock Accessible model categories.
	\newblock {\em Appl. Categ. Struct.}, 25(2):187--196, 2017.
	\newblock \href {https://doi.org/10.1007/s10485-015-9419-6}
	{\path{https://doi.org/10.1007/s10485-015-9419-6}}.
	
	\bibitem{rotho}
	J.~Rosick{\'{y}} and W.~Tholen.
	\newblock Left-determined model categories and universal homotopy theories.
	\newblock {\em {Trans. Am. Math. Soc.}}, 355(9):3611--3623, May 2003.
	\newblock \href {https://doi.org/10.1090/s0002-9947-03-03322-1}
	{\path{https://doi.org/10.1090/s0002-9947-03-03322-1}}.
	
	\bibitem{MR45:9323}
	R.~M. Vogt.
	\newblock Convenient categories of topological spaces for homotopy theory.
	\newblock {\em Arch. Math. (Basel)}, 22:545--555, 1971.
	\newblock \href {https://doi.org/10.1007/BF01222616}
	{\path{https://doi.org/10.1007/BF01222616}}.
	
	\bibitem{MR3722069}
	K.~Ziemia\'{n}ski.
	\newblock Spaces of directed paths on pre-cubical sets.
	\newblock {\em Appl. Algebra Engrg. Comm. Comput.}, 28(6):497--525, 2017.
	\newblock \href {https://doi.org/10.1007/s00200-017-0316-0}
	{\path{https://doi.org/10.1007/s00200-017-0316-0}}.
	
	\bibitem{MR4070250}
	K.~Ziemia\'{n}ski.
	\newblock Spaces of directed paths on pre-cubical sets {II}.
	\newblock {\em J. Appl. Comput. Topol.}, 4(1):45--78, 2020.
	\newblock \href {https://doi.org/10.1007/s41468-019-00040-z}
	{\path{https://doi.org/10.1007/s41468-019-00040-z}}.
	
\end{thebibliography}

\end{document}